\documentclass[11pt]{article}

\author{Claudia Kl{\"u}ppelberg\thanks{Center for Mathematical Sciences, Technical University of Munich, Boltzmannstra\ss e 3, 85748 Garching, Germany, e-mail: cklu@tum.de, vietson.pham@tum.de, URL: www.statistics.ma.tum.de}~~and Viet Son Pham$^\ast$}
\title{Estimation of causal CARMA random fields}
\date{\today}

\usepackage[T1]{fontenc}
\usepackage[utf8]{inputenc}
\usepackage{amsmath}
\usepackage{amsfonts}
\usepackage{amssymb}
\usepackage{amsthm}
\usepackage{amsopn}		
\usepackage{mathrsfs}	
\usepackage{dsfont}
\usepackage{color}
\usepackage{comment}
\usepackage[numbers]{natbib}
\usepackage{longtable}
\usepackage{enumitem}
\setenumerate{leftmargin=*}
\usepackage[linktocpage]{hyperref}
\usepackage{graphicx}
\usepackage{float}

\numberwithin{equation}{section}	
\allowdisplaybreaks[3]	
\newcommand{\bbn}{\mathbb{N}}
\newcommand{\bbz}{\mathbb{Z}}
\newcommand{\bbr}{\mathbb{R}}
\newcommand{\bbc}{\mathbb{C}}
\newcommand{\bbe}{\mathbb{E}}

\newcommand{\bbp}{\mathbb{P}}

\newcommand{\bone}{\mathds 1}

\newcommand{\calf}{{\cal F}}

\newcommand{\caln}{{\cal N}}
\newcommand{\calo}{{\cal O}}

\newcommand{\Del}{{\Delta}}
\newcommand{\ga}{{\gamma}}

\newcommand{\eps}{{\epsilon}}

\newcommand{\la}{{\lambda}}
\newcommand{\La}{{\Lambda}}
\newcommand{\si}{{\sigma}}

\newcommand{\om}{{\omega}}

\newcommand{\be}{\textbf{e}}
\newcommand{\bi}{\boldsymbol{i}}

\newcommand{\dd}{\mathrm{d}}
\newcommand{\ee}{\mathrm{e}}

\newcommand{\limd}{\stackrel{d}{\rightarrow}}	
\newcommand{\var}{{\mathrm{Var}}}	
\newcommand{\cov}{{\mathrm{Cov}}}	

\newcommand{\halmos}{\quad\hfill $\Box$}	

\newcommand{\Leb}{\mathrm{Leb}}

\newcommand{\argmin}{\mathrm{argmin}}
\renewcommand{\Re}{\mathrm{Re}}
\renewcommand{\Im}{\mathrm{Im}}

\newtheoremstyle{neu}
{11pt}      
{11pt}      
{}                  
{}          
{\bfseries} 
{}          
{1em}  
{\textbf{\thmname{#1}\thmnumber{ #2}\thmnote{ (#3)}}}          

\newtheoremstyle{proof}
{11pt}      
{11pt}      
{}                  
{}          
{\bfseries} 
{}            
{1em}          
{\textbf{\thmname{#1}.}}          

\newtheorem{Theorem}{Theorem}[section]

\newtheorem{Lemma}[Theorem]{Lemma}
\newtheorem{Proposition}[Theorem]{Proposition}
\theoremstyle{neu}
\newtheorem{Definition}[Theorem]{Definition}
\newtheorem{Example}[Theorem]{Example}
\newtheorem{Remark}[Theorem]{Remark}
\newtheorem{Algorithm}[Theorem]{Algorithm}
\newtheorem{Assumption}{Assumption}

\theoremstyle{proof}
\newtheorem{Proof}{Proof}

\begin{document}

\maketitle

\begin{abstract}
	We estimate model parameters of Lévy-driven causal CARMA random fields by fitting the empirical variogram to the theoretical counterpart using a weighted least squares (WLS) approach. Subsequent to deriving asymptotic results for the variogram estimator, we show strong consistency and asymptotic normality of the parameter estimator. Furthermore, we conduct a simulation study to assess the quality of the WLS estimator for finite samples. For the simulation we utilize numerical approximation schemes based on truncation and discretization of stochastic integrals and we analyze the associated simulation errors in detail. Finally, we apply our results to real data of the cosmic microwave background.
\end{abstract}

\vfill

\noindent
\begin{tabbing}
\end{tabbing}

\vspace{1cm}


\vspace{0.5cm}

\newpage

\section{Introduction}

Lévy-driven continuous-time autoregressive moving average (CARMA) processes are a well-studied class of stochastic processes and enjoy versatile applications in many disciplines (cf. \citet{Brockwell14} and the references therein). By contrast, considerably less is known about CARMA random fields indexed by $\bbr^d$, which have been defined only recently. To the best of our knowledge, two different classes exist in the literature: the \emph{isotropic CARMA random field} was introduced in \citet{Brockwell17} and the \emph{causal CARMA random field} in \cite{Pham18b}. While Bayesian parameter estimation is included in \cite{Brockwell17}, the paper \citet{Pham18b} only provides stochastic properties of causal CARMA random fields. The goal of this article is to provide a semiparametric method to estimate model parameters of causal CARMA random fields from discretely observed samples.

A Lévy-driven causal CARMA random field $(Y(t))_{t\in\bbr^d}$ on $\bbr^d$ is given by the equation
\begin{equation}\label{help5}
Y(t)=\int_{-\infty}^{t_1}\cdots\int_{-\infty}^{t_d} b^\top\ee^{A_1(t_1-s_1)}\cdots\ee^{A_d(t_d-s_d)}\be_p \,\La(\dd s),\quad t=(t_1,...,t_d)\in\bbr^d,
\end{equation}
where $A_1,...,A_d\in\bbr^{p\times p}$ are companion matrices, $\be_p=(0,...,0,1)^\top,b\in\bbr^p$ and $\La$ is a homogeneous Lévy basis, i.e., the multi-parameter analog of a Lévy process (see Section~\ref{sectionprelim} for more details). Due to its similar structure, many commonly known properties of CARMA processes also hold for $Y$, such as càdlàg sample paths, exponentially decreasing autocovariance functions and rational spectral densities. In fact, the random field $Y$ reduces to a causal CARMA process if $d=1$. Moreover, $Y$ has an autocovariance function which is both anisotropic and non-separable in the sense of \citet{Guttorp13}.

Since the matrices $A_1,...,A_d$ are in companion form, they are completely determined by their eigenvalues. These eigenvalues in conjunction with the components of the vector $b$ will form the model parameters. As our main tool for parameter estimation we choose the variogram, which is broadly applied in spatial statistics. It is defined as
\begin{equation*}
\psi(t)=\var[Y(t+s)-Y(s)],\quad t,s\in\bbr^d,
\end{equation*}
for stationary random fields (cf. Section~2.2.1 of \citet{Cressie93}). Furthermore, it is pointed out in Section~2.4.1 of \cite{Cressie93} that variogram estimation performs better than autocovariance estimation in terms of bias and in the presence of trend contamination. Assuming that observations of $Y$ are given on a regular lattice $L=\{ \Del,...,N\Del \}^d$, we estimate the model parameters by a two-step procedure.
First, we calculate an empirical version of the variogram $\psi(\cdot)$ at different lags using a non-parametric estimator $\psi_N^*(\cdot)$. Second, we fit the empirical variogram to the theoretical one using a weighted least squares method. More precisely, for a given set of strictly positive weights $w_j$, we estimate the true vector of CARMA parameters $\theta_0$ by means of the weighted least squares (WLS) estimator
\begin{equation*}
	\theta_N^*:=\argmin_{\theta\in\Theta}\Bigg\{ \sum_{j=1}^K w_j\left(\psi_N^*(t^{(j)})-\psi_\theta(t^{(j)})\right)^2 \Bigg\},
\end{equation*}
where $\Theta$ is a compact parameter space containing $\theta_0$ and $K$ is the number of lags used (see also Equation~\eqref{help13}).

An important task in connection with this approach is to determine sufficiently many lags $t^{(1)},...,t^{(K)}\in\bbr^d$ in order to obtain identifiability of the model parameters. We tackle this problem and show that under certain conditions a small number of lags on the principal axes of the Cartesian coordinate system is already sufficient to recover the CARMA parameters. In particular, one does not need to assume the property of \emph{invertibility} (or an analog thereof) as for CARMA processes. This fact differentiates the one-dimensional case from the higher dimensional case and we will investigate this in more detail.

Another part of this article is devoted to the study of different numerical simulation schemes for the causal CARMA random field. 
We derive approximation algorithms similar to those presented in \citet{Chen16} and \citet{Nguyen15} which are based on truncation or discretization of the stochastic integral in Equation~\eqref{help5}. We show that the output converges in mean-square and almost surely to the underlying CARMA random field. The algorithms are then used to conduct a simulation study in order to assess the quality of the WLS estimator. Subsequently, we apply the estimator to data of the cosmic microwave background.

Our paper is organized as follows: We recall the definition and basic properties of causal CARMA random fields in Section~\ref{sectionprelim}. Therein, a new formula for the spectral density is also proven. Strong consistency and asymptotic normality of the non-parametric variogram estimator $\psi_N^*(\cdot)$ is shown in Section~\ref{sectionempvario}. Subsequently, Section~\ref{sectionestimation} is concerned with the asymptotic properties of the WLS estimator $\theta_N^*$. Under identifiability conditions, we show strong consistency and asymptotic normality. While it is easier to show identifiability of CAR parameters, we obtain identifiability of CARMA parameters by carefully analyzing algebraic properties of the variogram. In Section~\ref{sectionsimu} we consider two different simulation methods and their associated algorithms. It is shown that the simulations converge pointwise both in $L^2$ and almost surely to the underlying true random fields as the truncation parameter tends to infinity and the discretization parameter tends to zero. The paper concludes with a simulation study and an application to cosmic microwave background data in Section~\ref{simstudy} and Section~\ref{sectionCMB}.

We use the following notation throughout this article: $\bone_{\{\cdot\}}$ denotes the indicator function such that for instance $\bone_{\{t\geq 0\}}$ is the Heaviside function. Furthermore, $A^\top$ denotes the transpose of a matrix (or a vector) $A$. The components of a vector $u\in\bbr^d$ are given by $u_1,...,u_d$ if not stated otherwise. For $u,v\in\bbr^d$, $\|u\|$ is the Euclidean norm, $u\cdot v\in\bbr$ is the scalar product, $u\odot v\in\bbr^d$ is the componentwise product, and $u\leq v$ if and only if $u_i\leq v_i$ for all $i\in\{1,...,d\}$. The $d$-dimensional interval $[u,v]$ is defined as $[u,v]:=\{s\in\bbr^d\colon u\leq s\leq v\}$ and we set $\bbr_+=[0,\infty)$. Additionally, $\be_1,...,\be_d$ are the unit vectors in $\bbr^d$ and $\be:=\be_1+\cdots+\be_d=(1,...,1)^\top$. Diagonal matrices are denoted by $\mathrm{diag}(\la_{1},...,\la_{d})\in\bbr^{d\times d}$ and $M_d(\bbr[z])$ is the space of all matrix polynomials of dimension $d\times d$. Finally, $\Re(z)$ and $\Im(z)$ are the real and imaginary part of a complex number $z$, $\Leb(\cdot)$ is the Lebesgue measure, and $\bi$ is the imaginary unit.

\section{Preliminaries}\label{sectionprelim}
First and foremost, we summarize some important properties of causal CARMA random fields. To this end, let us fix a probability space $(\Omega,\calf,\bbp)$, supporting all stochastic objects in this paper. As stated in the introduction, CARMA random fields are defined as stochastic integrals driven by \emph{homogeneous Lévy bases}. These are random measures which can be seen as a generalization of Lévy processes and their integration theory was developed in the seminal paper \citet{Rajput89}. For a homogeneous Lévy basis $\La$ we denote its characteristic triplet by $(\beta,\si^2,\nu)$, where $\beta\in\bbr$, $\si\in\bbr_+$ and $\nu$ is a Lévy measure. We say that $\La$ has a finite second moment and variance $\kappa_2:=\si^2 + \int_\bbr z^2\,\nu(\dd z)$ if and only if $\int_\bbr z^2\,\nu(\dd z)<\infty$. The variance $\kappa_2$ always appears in conjunction with the variogram or mean squared errors throughout this article. For more details on Lévy bases, we refer to Section~2 in \cite{Pham18b}. The following definition of causal CARMA random fields is taken from the same reference.
\begin{Definition}\label{defcarma} 
Let $q$ and $p$ be two non-negative integers such that $q<p$, $b=(b_0,...,b_{p-1})^\top\in\bbr^p$ with $b_q\neq 0$ and $b_i=0$ for $i>q$, $\be_p=(0,...,0,1)^\top\in\bbr^p$, and $A_i$ be the companion matrix to a monic polynomial $a_i$ of degree $p$ with real coefficients and roots having strictly negative real parts for $i=1,...,d$. A random field $(Y(t))_{t\in\bbr^d}$ is called \emph{(causal) CARMA$(p,q)$ random field} if it satisfies the equations
\begin{equation}\label{statespace}
\begin{aligned}
Y(t)&=b^\top X(t),\quad t\in\bbr^d,\\
X(t)&=\int_{-\infty}^{t_1}\cdots\int_{-\infty}^{t_d} \ee^{A_1(t_1-s_1)}\cdots\ee^{A_d(t_d-s_d)}\be_p \,\La(\dd s),\quad t\in\bbr^d,
\end{aligned}
\end{equation}
where $\La$ is a homogeneous Lévy basis on $\bbr^d$ with $\int_\bbr \log(|z|)^d\bone_{\{|z|>1\}} \,\nu(\dd z)<\infty$. A (causal) CARMA$(p,0)$ random field is also called a \emph{(causal) CAR$(p)$ random field}.
\halmos
\end{Definition}
Under the conditions specified in this definition, it was shown in \cite{Pham18b} that CARMA random fields exist and are well defined. Furthermore, they are by definition causal since the value of $Y(t)$ at $t\in\bbr^d$ only depends on the driving Lévy basis $\La$ on the set $(-\infty,t_1]\times\cdots\times(-\infty,t_d]$. This type of causality can be interpreted as a directional influence.
Also, it is immediate to see that they have the moving average representation
\begin{equation}
\label{convo} Y(t)=(g*\La)(t):=\int_{\bbr^d} g(t-s) \,\La(\dd s),\quad t\in\bbr^d,
\end{equation}
where the \emph{kernel} $g$ is given by
\begin{equation}\label{kernel} 
 g(s)=b^\top\ee^{A_1 s_1}\cdots\ee^{A_d s_d}\be_p\bone_{\{s\geq0\}},\quad s\in\bbr^d.
\end{equation}
The kernel $g$ is anisotropic in contrast to the isotropic CARMA random field in \cite{Brockwell17}. Additionally, it is non-separable, i.e., it cannot be written as a product of the form $g(s)=g_1(s_1)\cdots g_d(s_d)$ with real-valued functions $g_i$ except in the CAR$(1)$ case.
If $d=1$, we recover the classical kernel of a causal CARMA process (indexed by $\bbr$). 

\begin{Remark}
Causal CARMA random fields solve a system of stochastic partial differential equations, which generalizes the classical state-space representation of CARMA processes. For more details, see Section~3 in \cite{Pham18b}.
\halmos
\end{Remark}

In this article, we always impose the following additional conditions:
\begin{Assumption}~\label{assump}
\begin{itemize}
\item The Lévy basis $\La$ has mean zero and a finite second moment.
\item  The companion matrix $A_i$ has distinct eigenvalues for $i=1,...,d$.
\end{itemize}
\halmos
\end{Assumption}
The first part of Assumption~\ref{assump} ensures the existence of a second-order structure of $Y$, which is crucial for our estimation procedure in Section~\ref{sectionestimation}. In addition, the zero mean condition facilitates some computations, however, it is neither necessary nor restrictive. The autocovariance functions of CARMA random fields are non-separable (except in the CAR$(1)$ case), anisotropic, and integrable over $\bbr^d$ since they are exponentially decreasing (cf. \cite{Pham18b}). This implies the existence of a spectral density, which can be shown to be rational as for CARMA processes.

The second part of Assumption~\ref{assump} is analogous to Assumption~1 in \citet{Brockwell11}, where it is also pointed out that this condition is not critical since multiple eigenvalues can be handled as a limiting case.
Furthermore, this assumption implies that the kernel $g$ from Equation~\eqref{kernel} can alternatively be represented as
\begin{equation}\label{kernel2}
g(s)=\sum_{\la_1}\cdots\sum_{\la_d} d(\la_1,...,\la_d)\ee^{\la_1 s_1}\cdots \ee^{\la_d s_d}\bone_{\{s\geq0\}},\quad s\in\bbr^d,
\end{equation}
where $d(\la_1,...,\la_d)$ are (possibly complex) coefficients and $\sum_{\la_i}$ denotes the sum over distinct eigenvalues of $A_i$ for $i=1,...,d$ (cf. Corollary~3.7 in \cite{Pham18b}).

It is commonly known that an equidistantly sampled CARMA process is always an ARMA process.
Under certain conditions, we also obtain an ARMA random field if we sample a CARMA random field on a regular lattice (see Section~4.3 in \cite{Pham18b}). However, these conditions are rather restrictive (one of which is $A_1=\cdots=A_d$) and it is not known whether this sampling property can be generalized to the whole class of CARMA random fields. Nevertheless, we will see in Section~\ref{sectionsimu} that a CARMA random field sampled on a regular grid can always be approximated arbitrarily well by a discrete-parameter moving average random field (of finite order) in terms of the mean squared error and almost surely.

From Equation~\eqref{convo} we observe that $Y$ is strictly stationary, which in turn implies that the variogram
\begin{equation*}
\psi(t)=\var[Y(t)-Y(0)]=\var[Y(t+s)-Y(s)],\quad t,s\in\bbr^d,
\end{equation*}
is translation-invariant, i.e., independent of $s$.
In Section~\ref{sectionestimation} we will estimate the CARMA parameters $b$ and the eigenvalues of $A_1,...,A_d$ by fitting an empirical version of $\psi$ to its theoretical counterpart. Therefore, it is necessary to have the variogram structure of CARMA random fields at hand, which is given in the next proposition.
\begin{Proposition}\label{vario}
Suppose that $(Y(t))_{t\in\bbr^d}$ is a CARMA$(p,q)$ random field such that Assumption~\ref{assump} holds true. Then the variogram $\psi$ of $Y$ has the form
\begin{align*}
\psi(t)&=2\kappa_2\sum_{\la_1}\cdots\sum_{\la_d}\sum_{v\in\{-1,1\}^d} d_v(\la_1,...,\la_d)\bone_{\{ t\odot v\in\bbr_+^d \}} \left(1-\ee^{\la_1|t_1|}\cdots \ee^{\la_d|t_d|}\right),\quad t\in\bbr^d,
\end{align*}
where $\{d_v(\la_1,...,\la_d)\}$ is a set of complex coefficients for every $v\in\{-1,1\}^d$ such that
\begin{equation*}
d_v(\la_1,...,\la_d)=d_{-v}(\la_1,...,\la_d)
\end{equation*}
and $\sum_{\la_i}$ denotes the sum over distinct eigenvalues of $A_i$ for $i=1,...,d$.
\end{Proposition}
\begin{Proof}
The statement is a combination of Theorem~4.1. in \cite{Pham18b} and the relation
\begin{equation*}
\psi(t)=2(\ga(0)-\ga(t)),\quad t\in\bbr^d,
\end{equation*}
where
\begin{equation}\label{help9}
\ga(t)=\cov[Y(t),Y(0)]=\cov[Y(t+s),Y(s)],\quad t,s\in\bbr^d,
\end{equation}
is the autocovariance function of $Y$.
\halmos
\end{Proof}
As it was argued in \cite{Pham18b}, it is in general hard to find explicit formulae for $d_v(\la_1,...,\la_d)$ in terms of $d,p,q,b,A_1,...,A_d$. However, if we fix the dimension $d$ and the orders $p$ and $q$, we are able to compute the variogram explicitly. We consider the following example.
\begin{Example}
Let $d=2$, $p=2$, $q=1$ and $\La$ be a homogeneous Lévy basis satisfying Assumption~\ref{assump}.
We assume that the CARMA$(2,1)$ random field
\begin{equation*}
Y(t)=\int_{-\infty}^{t_1}\int_{-\infty}^{t_2} b^\top\ee^{A_1(t_1-s_1)}\ee^{A_2(t_2-s_2)}\be_p \,\La(\dd s),\quad t\in\bbr^2,
\end{equation*}
has parameters $b=(b_0,b_1)\in\bbr^2$,
\begin{equation*}
A_1=\begin{pmatrix}0&1\\-\la_{11}\la_{12}&\la_{11}+\la_{12}\end{pmatrix}\quad\text{and}\quad A_2=\begin{pmatrix}0&1\\-\la_{21}\la_{22}&\la_{21}+\la_{22}\end{pmatrix},
\end{equation*}
such that the eigenvalues $\la_{11},\la_{12},\la_{21},\la_{22}\in\bbr$ have strictly negative real parts and satisfy $\la_{11}\neq\la_{12}$ and $\la_{21}\neq\la_{22}$. In this case, the variogram $\psi$ of $Y$ is given by
\begin{equation*}
\psi(t)=2\kappa_2\sum_{k=11}^{12}\sum_{l=21}^{22} \Big( d_{(1,1)}(\la_k,\la_l)\bone_{\{ t_1t_2\geq0 \}}
+d_{(1,-1)}(\la_k,\la_l)\bone_{\{ t_1t_2<0 \}} \Big) \left( 1-\ee^{\la_k|t_1|}\ee^{\la_l|t_2|} \right),\quad t\in\bbr^2,
\end{equation*}
where
\begin{align*}
d_{(1,1)}(\la_{11},\la_{21})&=\frac{(\la_{12}-\la_{21}) (b_0+b_1 \la_{11}) (b_0 (2 \la_{11}+\la_{12}+\la_{21})+b_1 \la_{11} (\la_{12}-\la_{21}))}{4 \la_{11}
   \la_{21} (\la_{11}-\la_{12}) (\la_{11}+\la_{12}) (\la_{21}-\la_{22}) (\la_{21}+\la_{22})}\\  
d_{(1,1)}(\la_{12},\la_{21})&=-\frac{(\la_{11}-\la_{21}) (b_0+b_1 \la_{12}) (b_0 (\la_{11}+2 \la_{12}+\la_{21})+b_1 \la_{12} (\la_{11}-\la_{21}))}{4 \la_{12}
   \la_{21} (\la_{11}-\la_{12}) (\la_{11}+\la_{12}) (\la_{21}-\la_{22}) (\la_{21}+\la_{22})}\\  
d_{(1,1)}(\la_{11},\la_{22})&=\frac{(\la_{12}-\la_{22}) (b_0+b_1 \la_{11}) (b_0 (2 \la_{11}+\la_{12}+\la_{22})+b_1 \la_{11} (\la_{12}-\la_{22}))}{4 \la_{11}
   \la_{22} (\la_{11}-\la_{12}) (\la_{11}+\la_{12}) (\la_{22}-\la_{21}) (\la_{21}+\la_{22})}\\
d_{(1,1)}(\la_{12},\la_{22})&=\frac{(\la_{11}-\la_{22}) (b_0+b_1 \la_{12}) (b_0 (\la_{11}+2 \la_{12}+\la_{22})+b_1 \la_{12} (\la_{11}-\la_{22}))}{4 \la_{12}
   \la_{22} (\la_{11}-\la_{12}) (\la_{11}+\la_{12}) (\la_{21}-\la_{22}) (\la_{21}+\la_{22})}
\end{align*}
and
\begin{align*}
d_{(1,-1)}(\la_{11},\la_{21})&=\frac{(\la_{12}+\la_{21}) (b_0+b_1 \la_{11}) (b_0 (2 \la_{11}+\la_{12}-\la_{21})+b_1 \la_{11} (\la_{12}+\la_{21}))}{4 \la_{11}
   \la_{21} (\la_{11}-\la_{12}) (\la_{11}+\la_{12}) (\la_{21}-\la_{22}) (\la_{21}+\la_{22})}\\
d_{(1,-1)}(\la_{12},\la_{21})&=-\frac{(\la_{11}+\la_{21}) (b_0+b_1 \la_{12}) (b_0 (\la_{11}+2 \la_{12}-\la_{21})+b_1 \la_{12} (\la_{11}+\la_{21}))}{4 \la_{12}
   \la_{21} (\la_{11}-\la_{12}) (\la_{11}+\la_{12}) (\la_{21}-\la_{22}) (\la_{21}+\la_{22})}\\
d_{(1,-1)}(\la_{11},\la_{22})&=-\frac{(\la_{12}+\la_{22}) (b_0+b_1 \la_{11}) (b_0 (2 \la_{11}+\la_{12}-\la_{22})+b_1 \la_{11} (\la_{12}+\la_{22}))}{4 \la_{11}
   \la_{22} (\la_{11}-\la_{12}) (\la_{11}+\la_{12}) (\la_{21}-\la_{22}) (\la_{21}+\la_{22})}\\
d_{(1,-1)}(\la_{12},\la_{22})&=\frac{(\la_{11}+\la_{22}) (b_0+b_1 \la_{12}) (b_0 (\la_{11}+2 \la_{12}-\la_{22})+b_1 \la_{12} (\la_{11}+\la_{22}))}{4 \la_{12}
   \la_{22} (\la_{11}-\la_{12}) (\la_{11}+\la_{12}) (\la_{21}-\la_{22}) (\la_{21}+\la_{22})}
\end{align*}
These formulae have been computed with the computer algebra system \texttt{Mathematica}.
\halmos
\end{Example}
The next result, which is of theoretical interest and will be useful later on, contains a formula for the spectral density of $Y$ which is more explicit than Equation~(4.5) in \cite{Pham18b}.
\begin{Proposition}\label{specden}
Suppose that $(Y(t))_{t\in\bbr^d}$ is a CARMA$(p,q)$ random field such that $\La$ has a finite second moment. Further, let $a_i(z)=\sum_{j=0}^{p} a_{i,j} z^{p-j}$ for $i=1,...,d$ be the monic polynomials in Definition~\ref{defcarma}. Then the spectral density $f$ of $Y$ has the representation
\begin{equation*}
f(\om)=\frac{\kappa_2 }{(2\pi)^d }\left| \frac{Q(\bi\om)}{P(\bi\om)}\right|^2,\quad \om\in\bbr^d,
\end{equation*}
with polynomials
\begin{equation*}
P(z)=a_1(z_1)\cdots a_d(z_d),
\end{equation*}
and
\begin{equation*}
Q(z)=b^\top Q_1(z_1)\cdots Q_d(z_d)c,
\end{equation*}
and matrix polynomials
\begin{equation*}
Q_i(z)=a_i(z)(zI_p-A_i)^{-1}\in M_p(\bbr[z]).
\end{equation*}
For each $i=1,...,d$, the $(k,l)$-entry of the matrix polynomial $Q_i$ is given by
\begin{equation*}
Q_{i,k,l}(z)=\begin{cases}
z^{p-1+k-l} + \sum_{j=1}^{p-l}a_{i,j}z^{p-1-j+k-l} &k\leq l,\\
-\sum_{j=p-l+1}^{p}a_{i,j}z^{p-1-j+k-l} &k>l.
\end{cases}
\end{equation*}
\end{Proposition}
\begin{Proof}
Proposition~11.2.2 in \citet{Bernstein05} and Equation~\eqref{kernel} imply that the Fourier transform of the kernel $g$ satisfies
\begin{equation*}
\tilde g(\om)=\int_{\bbr^d} g(s)\ee^{-\bi \om\cdot s}\,\dd s=b^\top(\bi\om_1I_p-A_1)^{-1}\cdots(\bi\om_dI_p-A_d)^{-1}\be_p,\quad \om\in\bbr^d.
\end{equation*}
Applying Lemma~3.1 in \citet{Brockwell13b} and the relation
\begin{equation*}
f(\om)=\frac{\kappa_2 }{(2\pi)^d }\tilde g(\om)\tilde g(-\om)=\frac{\kappa_2 }{(2\pi)^d }|\tilde g(\om)|^2,
\end{equation*}
yields the claimed assertion.
\halmos
\end{Proof}

\section{Asymptotic properties of the empirical variogram}\label{sectionempvario}

Let  $(Y(t))_{t\in\bbr^d}$ be a CARMA$(p,q)$ random field satisfying Assumption~\ref{assump}. If we are given observations of $Y$ on a lattice $L=\{ \Del,...,N\Del \}^d$, we can estimate the variogram $\psi(\cdot)$ by Matheron's method-of-moment estimator (cf. Section~2.4 in \citet{Cressie93} for more details)
\begin{equation*}
\psi_N^*(t):=\frac{1}{|B_{N,t}|}\sum_{s\in B_{N,t}}(Y(t+s)-Y(s))^2, \quad t\in\{ (1-N)\Del,...,(N-1)\Del \}^d,
\end{equation*}
where
\begin{equation*}
B_{N,t}:=\{ s\in\Del\bbz^d\colon s,s+t\in\{ \Del,...,N\Del \}^d \}\text{ and }|B_{N,t}|=\prod_{i=1}^d (N-|t_i|)\bone_{\{|t_i|\leq N\}}.
\end{equation*}
We aim to show strong consistency and multivariate asymptotic normality of $\psi_N^*(\cdot)$ as $N$ tends to infinity. To this end, we make use of the asymptotic normality of the autocovariance estimator
\begin{equation*}
\ga_N^*(t):=\frac{1}{|B_{N,t}|}\sum_{s\in B_{N,t}}Y(t+s)Y(s),
\end{equation*}
which was shown in \cite{Berger17}
for moving average random fields by applying a blocking technique and a central limit theorem for \emph{m-dependent} random fields.
\begin{Theorem}\label{thmvario}
Suppose that $(Y(t))_{t\in\bbr^d}$ is a CARMA$(p,q)$ random field such that Assumption~\ref{assump} holds true, $\La$ has a finite fourth moment $\kappa_4$ and observations of $Y$ are given on the lattice $L=\{ \Del,...,N\Del \}^d$. Then we have for all $t\in\Del\bbz^d$ that
\begin{equation*}
\lim_{N\to\infty}\psi_N^*(t)=\psi(t) \quad\text{a.s.}.
\end{equation*}
Further let $t^{(1)},...,t^{(K)}\in\Del\bbz^d$ be $K$ distinct lags and $t^{(0)}=(0,...,0)^\top$. Then we have
\begin{equation*}
N^{d/2}(\psi_N^*(t^{(1)})-\psi(t^{(1)}),...,\psi_N^*(t^{(K)})-\psi(t^{(K)}))\limd\caln(0,FVF^\top),\quad\text{as } N\to\infty,
\end{equation*}
where the two matrices $F$ and $V=(v_{i,j})_{i,j=0,...,K}$ are given by
\begin{equation*}
F=2\begin{pmatrix}
1		&-1		&0 	&\cdots	&0\\
1		&0		&-1 &\cdots	&0\\
\vdots  &\vdots &\vdots	&\ddots	&\vdots\\
1		&0		&0 	&\cdots	&-1
\end{pmatrix}\in\bbr^{K\times(K+1)}
\end{equation*}
and
\begin{equation*}
v_{i,j}=\sum_{l\in\Del\bbz^d}(\kappa_4-3\kappa_2^2)\int_\bbr^d g(s)g(s+t^{(i)})g(s+l)g(s+l+t^{(j)})\,\dd s +\ga(l)\ga(l+t^{(i)}-t^{(j)}) +\ga(l+t^{(i)})\ga(l-t^{(j)})
\end{equation*}
for all $i,j=0,...,K$.
\end{Theorem}
\begin{Proof}
First of all, we show strong consistency of the variogram estimator $\psi_N^*(\cdot)$. By Corollary~3.18 in \cite{Berger17}, we have for all $t\in\Del\bbz^d$ that
\begin{equation*}
\lim_{N\to\infty}\ga_N^*(t)=\ga(t)\quad\text{a.s.}.
\end{equation*}
Considering the following limit
\begin{align*}
&\lim_{N\to\infty}\psi_N^*(t)-2(\ga_N^*(0)-\ga_N^*(t))\\
&\quad=\lim_{N\to\infty}\left(\frac{1}{|B_{N,t}|}\sum_{s\in B_{N,t}}(Y(t+s)^2+Y(s)^2) - \frac{1}{|B_{N,0}|}\sum_{s\in B_{N,0}}Y(s)^2\right)\\
&\quad= 2\ga(0)-2\ga(0)=0\quad\text{a.s.},
\end{align*}
we deduce that
\begin{equation*}
\lim_{N\to\infty}\psi_N^*(t)=\lim_{N\to\infty}2(\ga_N^*(0)-\ga_N^*(t))=2(\ga(0)-\ga(t))=\psi(t) \quad\text{a.s.}
\end{equation*}
as desired. It remains to show asymptotic normality of $\psi_N^*(\cdot)$.
Since the kernel $g$ in Equation~\eqref{kernel2} is a sum of exponentials, we have for every $i,j=0,...,K $ that
\begin{equation*}
\int_{\bbr^d}|g(s)g(s+k)g(s+t^{(i)})g(s+k+t^{(j)})|\,\dd s = \calo(\ee^{2\la_{max,1}k_1+...+2\la_{max,1}k_d}),\quad k\to\infty,
\end{equation*}
where $\la_{max,i}:=\max\{ \Re(\la_i) \colon \la_i \text{ is eigenvalue of } A_i \}<0$ for $i=1,...,d$. Hence, we obtain
\begin{equation*}
\sum_{k\in\Del\bbz^d}\int_{\bbr^d}|g(s)g(s+k)g(s+t^{(i)})g(s+k+t^{(j)})|\,\dd s<\infty.
\end{equation*}
Moreover, by Theorem~4.1. in \cite{Pham18b} we also have that
\begin{equation*}
\sum_{k\in\Del\bbz^d}\ga(k)^2<\infty.
\end{equation*}
We conclude that the conditions of Theorem~3.8 in \cite{Berger19} are satisfied, which in turn shows that
\begin{equation}\label{help8}
N^{d/2}(\ga_N^*(t^{(0)})-\ga(t^{(0)}),...,\ga_N^*(t^{(K)})-\ga(t^{(K)}))\limd\caln(0,V),\quad\text{as } N\to\infty.
\end{equation}
Consider now the mapping
\begin{equation*}
f\colon\bbr^{K+1}\to\bbr^K, \quad f(x_0,...,x_K)\mapsto(2(x_0-x_1),...,2(x_0-x_K)),
\end{equation*}
whose Jacobian is the matrix $F$.
The multivariate delta method (see e.g. Proposition~6.4.3 in \cite{Brockwell91}) in combination with the mapping $f$ and \eqref{help8} yields
\begin{equation*}
N^{d/2}(2(\ga_N^*(0)-\ga_N^*(t^{(1)}))-\psi(t^{(1)}),...,2(\ga_N^*(0)-\ga_N^*(t^{(K)}))-\psi(t^{(K)}))\limd\caln(0,FVF^\top)
\end{equation*}
as $N\to\infty$, and Slutsky's theorem finishes the proof.
\halmos
\end{Proof}

\section{Estimation of CARMA random fields}\label{sectionestimation}

According to Definition~\ref{defcarma}, a CARMA random field is determined by the pair $(p,q)$, the vector $b$, the companion matrices $A_1,...,A_d$ and the Lévy basis $\La$. To avoid redundancies in model specification one usually assumes that either $b_0$ or $\kappa_2$ is known. We assume the latter and thus the goal of this section is to estimate $b$ and $A_1,...,A_d$ when $p$, $q$ and $\kappa_2$ are given. Since every companion matrix is uniquely determined by its eigenvalues, we define the \emph{CARMA parameter vector} $\theta$ as
\begin{equation}\label{help13}
\theta=(b_0,...,b_q,\la_{11},...,\la_{1p},\la_{21},...,\la_{dp})\in\bbr^{q+1}\times\bbc^{dp},
\end{equation}
where $\la_{i1},...,\la_{ip}$ are the eigenvalues of $A_i$ for $i=1,...,d$. Recall that $A_i$ is real by definition and thus its eigenvalues are real or appear in pairs of complex conjugates.
In order to estimate $\theta$, we fit the empirical variogram $\psi_N^*(\cdot)$ of the last section to the theoretical variogram $\psi_\theta(\cdot)$ using a weighted least squares approach. In other words, we consider the estimator
\begin{equation}\label{WLSest}
\theta_N^*:=\argmin_{\theta\in\Theta}\Bigg\{ \sum_{j=1}^K w_j\left(\psi_N^*(t^{(j)})-\psi_\theta(t^{(j)})\right)^2 \Bigg\}
\end{equation}
where $\Theta\subseteq\bbr^{q+1}\times\bbc^{dp}$ is a compact parameter space containing the true parameter vector $\theta_0$, $w_j>0$ are strictly positive weights and $t^{(1)},...,t^{(K)}\in\bbr^d$ are prescribed lags. The paper \citet{Lahiri02} determines asymptotic properties of least squares estimators for parametric variogram models subject to asymptotic properties of the underlying variogram estimators. We use these results in conjunction with Theorem~\ref{thmvario} to show strong consistency and asymptotic normality of $\theta_N^*$. In the following, we denote by
\begin{equation*}
\xi_i(\theta):=((\partial/\partial\theta_i)\psi_\theta(t^{(1)}),...,(\partial/\partial\theta_i)\psi_\theta(t^{(K)}))
\end{equation*}
the vector of first order partial derivatives of $\psi_\theta(t^{(1)}),...,\psi_\theta(t^{(K)})$ with respect to the $i$'th coordinate of $\theta$ and define
\begin{equation*}
\Xi(\theta):=-(\xi_1(\theta),...,\xi_{dp+q+1}(\theta)),
\end{equation*}
which is the Jacobian matrix of the mapping $\theta\mapsto(\psi_\theta(t^{(1)}),...,\psi_\theta(t^{(K)}))$.
\begin{Theorem}\label{thmWLS}
Suppose that $(Y(t))_{t\in\bbr^d}$ is a CARMA$(p,q)$ random field with true parameter vector $\theta_0$ such that Assumption~\ref{assump} holds true, $\kappa_2=1$, $\La$ has a finite fourth moment and observations of $Y$ are given on the lattice $L=\{ \Del,...,N\Del \}^d$. Further, assume that
\begin{itemize}
\item the true parameter vector $\theta_0$ lies inside a compact parameter space $\Theta\subseteq\bbr^{q+1}\times\bbc^{dp}$,
\item the mapping $\Theta\ni\theta\mapsto (\psi_\theta(t^{(1)}),...,\psi_\theta(t^{(K)}))$ is injective (identifiability criterion).
\end{itemize}
Then we have both
\begin{equation*}
\lim_{N\to\infty}\theta_N^*=\theta_0 \quad\text{a.s.}
\end{equation*}
and
\begin{equation*}
N^{d/2}(\theta_N^*-\theta_0)\limd\caln(0,\Sigma),\quad\text{as } N\to\infty,
\end{equation*}
where
\begin{equation*}
\Sigma=B(\theta_0)\Xi(\theta_0)^\top W FVF^\top W\Xi(\theta_0)B(\theta_0),
\end{equation*}
with $F$ and $V$ as in Theorem~\ref{thmvario}, $W=diag(w_1,...,w_K)$ and $B(\theta_0)=(\Xi(\theta_0)^\top W \Xi(\theta_0))^{-1}$.
\end{Theorem}
\begin{Proof}
We only have to check conditions (C.1)-(C.3) in \cite{Lahiri02} since our assertions follow directly from Theorems~3.1 and 3.2 of this reference. 
Since $\psi_\theta(t)=2(\ga_\theta(0)-\ga_\theta(t))$, it suffices to show that for each $t\in\bbr^d$ the autocovariance $\ga_\theta(t)$ is continuously differentiable with respect to $\theta$ in order to check (C.2)(ii). Recall that the relation
\begin{equation*}
A_i=\begin{pmatrix}
0	&1	&0	&\cdots	&0\\
0	&0	&1	&\cdots	&0\\
\vdots	&\vdots	&\vdots	&\ddots	&\vdots\\
0	&0	&0	&\cdots	&1\\
-a_{ip}	&-a_{i(p-1)}	&-a_{i(p-2)}	&\cdots	&-a_{i1}
\end{pmatrix}=V\mathrm{diag}(\la_{i1},...,\la_{ip})V^{-1}
\end{equation*}
is satisfied for companion matrices, where $V$ is the Vandermonde matrix
\begin{equation*}
V=\begin{pmatrix}
1	&\cdots	&1\\
\la_{i1}	&\cdots	&\la_{ip}\\
\vdots	&\ddots	&\vdots\\
\la_{i1}^{p-1}	&\cdots	&\la_{ip}^{p-1}
\end{pmatrix}.
\end{equation*}
Now assume $t\in\bbr_+^d$ first. Then (the proof of) Theorem~4.1. in \cite{Pham18b} implies
\begin{equation}\label{help10}
\ga_\theta(t)=\kappa_2b^\top \left( \int_{\bbr_+^2} \ee^{A_1 (s_1+t_1)}\cdots\ee^{A_d (s_d+t_d)}\be_p\be_p^\top\ee^{A^\top_d s_d}\cdots\ee^{A^\top_1 s_1} \,\dd s \right) \,b.
\end{equation}
Owing to the exponential structure of the integrand we recognize that $\ga_\theta(t)$ is in fact infinitely often differentiable with respect to $\theta$ and therefore condition (C.2)(ii) holds true for every $t\in\bbr_+^d$. For each $t\in\bbr^d$ an analogous argument applies with a slightly different integrand. Moreover, condition (C.2)(ii) implies both (C.2)(i) and (C.1) in light of the identifiability criterion and the fact that $\Theta$ is compact. Finally, the condition (C.3) is trivial since $W$ does not depend on $\theta$.
\halmos
\end{Proof}
An important task in connection with the previous theorem is to determine a sufficient set of lags $t^{(1)},...,t^{(K)}$ such that the identifiability condition is satisfied. For CAR$(p)$ random fields it is enough to consider finitely many lags on the principal axes of the Cartesian coordinate system. Before examining this matter, we prepend an auxiliary lemma which presents a simplified representation of the variogram on the principal axes.
\begin{Lemma}\label{lem1}
Suppose that $(Y(t))_{t\in\bbr^d}$ is a CARMA$(p,q)$ random field such that Assumption~\ref{assump} holds true. Then there exists a set of complex coefficients $\{ d^*_i(\la_i)\colon \la_i\text{ is an eigenvector of }A_i, i=1,...,d \}$ such that the values of the variogram $\psi$ is given on the principal axes by 
\begin{equation}\label{help12}
\psi(\tau\be_i )=2\kappa_2\sum_{\la_i}d^*_i(\la_i) \left(1-\ee^{\la_i|\tau|}\right),\quad \tau\in\bbr,\quad i=1,...,d,
\end{equation}
where $\sum_{\la_i}$ denotes the sum over distinct eigenvalues of $A_i$.
\end{Lemma}
\begin{Proof}
Proposition~\ref{vario} implies for every $i=1,...,d$ and $\tau\in\bbr^d$ that
\begin{align*}
\psi(\tau\be_i)&=2\kappa_2\sum_{\la_1}\cdots\sum_{\la_d} d_\be(\la_1,...,\la_d) \left(1-\ee^{0}\cdots \ee^{\la_i|\tau|}\cdots\ee^{0}\right)\\
&=2\kappa_2\sum_{\la_i}\left( \sum_{\la_j\colon j=1,...,i-1,i+1,...,d} d_\be(\la_1,...,\la_d) \right)\left(1-\ee^{\la_i|\tau|}\right)\\
&=2\kappa_2\sum_{\la_i}d^*_i(\la_i) \left(1-\ee^{\la_i|\tau|}\right).
\end{align*}
\halmos
\end{Proof}
The next example displays more explicit formulae for $d^*_i(\la_i)$ in the case of a CARMA$(2,1)$ random field.
\begin{Example}\label{varioaxes}
Let $(Y(t))_{t\in\bbr^2}$ be a CARMA$(2,1)$ random field such that Assumption~\ref{assump} holds true. Then we have for every $t_1,t_2\in\bbr$ that
\begin{align*}
\psi(t_1,0)&=\frac{ 2\kappa_2(b_0+b_1 \la_{12}) \left(b_0 \left(\la_{11}^2+2 \la_{11} \la_{12}+\la_{21} \la_{22}\right)+b_1 \la_{12}
\left(\la_{11}^2-\la_{21} \la_{22}\right)\right)}{4 \la_{12} \la_{21} \la_{22} (\la_{11}-\la_{12}) (\la_{11}+\la_{12})
(\la_{21}+\la_{22})}(1-e^{\la_{12} |t_1|})\\
&\quad+\frac{ 2\kappa_2(b_0+b_1 \la_{11}) \left(b_0 \left(2 \la_{11} \la_{12}+\la_{12}^2+\la_{21}
\la_{22}\right)+b_1 \la_{11} \left(\la_{12}^2-\la_{21} \la_{22}\right)\right)}{4 \la_{11} \la_{21} \la_{22} (\la_{12}-\la_{11})
(\la_{11}+\la_{12}) (\la_{21}+\la_{22})}(1-e^{\la_{11} |t_1|})
\end{align*}
and
\begin{align*}
\psi(0,t_2)&=2\kappa_2\frac{ b_0^2 \left(\la_{11}^2+3 \la_{11} \la_{12}+\la_{12}^2-\la_{21}^2\right)+2 b_0 b_1 \la_{11} \la_{12}
(\la_{11}+\la_{12})+b_1^2 \la_{11} \la_{12} \left(\la_{11} \la_{12}-\la_{21}^2\right)}{4 \la_{11} \la_{12} \la_{21}
(\la_{11}+\la_{12}) (\la_{22}-\la_{21}) (\la_{21}+\la_{22})}\\
&\quad \times(1-e^{\la_{21} |t_2|})\\
&\quad+2\kappa_2\frac{ b_0^2 \left(\la_{11}^2+3 \la_{11}
\la_{12}+\la_{12}^2-\la_{22}^2\right)+2 b_0 b_1 \la_{11} \la_{12} (\la_{11}+\la_{12})+b_1^2 \la_{11} \la_{12} \left(\la_{11}
\la_{12}-\la_{22}^2\right)}{4 \la_{11} \la_{12} \la_{22} (\la_{11}+\la_{12}) (\la_{21}-\la_{22}) (\la_{21}+\la_{22})}\\
&\quad \times(1-e^{\la_{22} |t_2|}).
\end{align*}
\halmos
\end{Example}
The next theorem establishes the identifiability of CAR$(p)$ parameters. Note that replacing the vector $b$ by $-b$ would not change the variogram. Hence, we may assume that $b_0$ is non-negative.
\begin{Theorem}\label{identCAR}
Suppose that $(Y(t))_{t\in\bbr^d}$ is a CAR$(p)$ random field such that Assumption~\ref{assump} holds true, $\kappa_2$ is given, $b_0\geq 0$, all eigenvalues $\la$ of $A_1,...,A_d$ satisfy $-\pi/\Del\leq\Im(\la)<\pi/\Del$ and all coefficients $d^*_i(\la_i)$ in Lemma~\ref{lem1} are nonzero. Then $\theta$ is uniquely determined by the variogram ordinates $\{ \psi(j\Del\be_i)\colon i=1,...,d;j=0,...,2p+1 \}$.
\end{Theorem}
\begin{Proof}
Assuming without loss of generality that $\Del=1$ and $\kappa_2=1/2$, and setting $\la_{i0}=0$ and $d_i^*(\la_{i0})=-\sum_{k=1}^p d_i^*(\la_{ik})$, Lemma~\ref{lem1} implies that 
\begin{equation*}
\psi(j\be_i)=-\sum_{k=0}^p d_i^*(\la_{ik}) \ee^{\la_{ik}j},\quad i=1,...,d,\quad j=0,...,2p+1.
\end{equation*}
Note that $-d_i^*(\la_{i0})=\psi(0)$ is twice the variance of $Y$ and therefore nonzero. Introducing the polynomials
\begin{equation*}
R_i(z):=\prod_{l=0}^p (z-\ee^{\la_{il}})=:\sum_{l=0}^{p+1} r_{il}z^l,
\end{equation*}
we observe for each $i=1,...,d$ and $j=0,...,2p+1$ that
\begin{equation*}
\sum_{l=0}^{p+1} r_{il}\psi((j+l)\be_i)=-\sum_{l=0}^{p+1} r_{il}\sum_{k=0}^p d_i^*(\la_{ik}) \ee^{\la_{ik}j}\ee^{\la_{ik}l}=-\sum_{k=0}^p d_i^*(\la_{ik}) \ee^{\la_{ik}j} \sum_{l=0}^{p+1} r_{il} \ee^{\la_{ik}l}=0,
\end{equation*}
where the last equation follows from the definition of $R_i(z)$. Hence, we get the linear systems
\begin{equation}\label{help15}
\Psi_i \begin{pmatrix}r_{i0}\\ \vdots \\ r_{ip}\end{pmatrix}:=
\begin{pmatrix}
\psi(0\be_i)	&\cdots	&\psi(p\be_i)\\
\vdots	&\ddots	&\vdots\\
\psi(p\be_i)	&\cdots	&\psi(2p\be_i)
\end{pmatrix}
\begin{pmatrix}r_{i0}\\ \vdots \\ r_{ip}\end{pmatrix}
=-\begin{pmatrix}\psi((p+1)\be_i)\\ \vdots \\ \psi((2p+1)\be_i)\end{pmatrix},
\end{equation}
where the system matrices $\Psi_i$ are quadratic Hankel matrices. We show that all $\Psi_i$ are invertible. To this end, 
for fixed $i\in\{1,...,d\}$,
assume that there is a vector $u=(u_0,...,u_p)\in\bbr^{p+1}$ satisfying
\begin{equation*}
\Psi_i u=0,
\end{equation*}
that is $u$ is an element inside $\Psi_i$'s kernel. Defining the polynomial $P(z):=\sum_{l=0}^{p} u_l z^l$, we obtain for all $j=0,...,p$ that
\begin{align*}
0=\sum_{l=0}^{p} u_l \psi((l+j)\be_i)
=\sum_{l=0}^{p} u_l \left( -\sum_{k=0}^p d_i^*(\la_{ik}) \ee^{\la_{ik}(l+j)} \right) = -\sum_{k=0}^p \ee^{\la_{ik}j} d_i^*(\la_{ik}) P(\ee^{\la_{ik}}).
\end{align*}
This gives the linear system
\begin{equation*}
\begin{pmatrix}
1	&\cdots	&1\\
\ee^{\la_{i0}}	&\cdots	&\ee^{\la_{ip}}\\
\vdots	&\ddots	&\vdots\\
\ee^{\la_{i0}p}	&\cdots	&\ee^{\la_{ip}p}
\end{pmatrix}
\begin{pmatrix} d_i^*(\la_{i0}) P(\ee^{\la_{i0}}) \\ d_i^*(\la_{i1}) P(\ee^{\la_{i1}}) \\ \vdots \\ d_i^*(\la_{ip}) P(\ee^{\la_{ip}})\end{pmatrix}
=0.
\end{equation*}
Since the system matrix is a regular Vandermonde matrix and the coefficients $d_i^*(\la_{ik})$ are nonzero, we conclude that  $P(\ee^{\la_{ik}})=0$ for $k=0,...,p$. The polynomial $P(\cdot)$ has $(p+1)$ different roots and is of degree $p$. Consequently, it has to be the zero polynomial, which means that $u=0$ and $\Psi_i$ is invertible.
By solving the linear systems \eqref{help15} we get all $r_{il}$, which gives the $\ee^{\la_{il}}$ by determining the roots of $R_i(z)$. Finally, all eigenvalues $\la_{il}$ can be obtained uniquely using the condition on the imaginary part $\Im(\la)$, and it is trivial to recover $b_0$ in light of Equations~\eqref{help9} and \eqref{help10}.
\halmos
\end{Proof}

\begin{Remark}
\begin{enumerate}
\item The set of parameter vectors $\theta$ of CAR random fields which have at least one vanishing coefficient $d^*_i(\la_i)$ is a lower dimensional algebraic variety in the parameter space $\bbr\times\bbc^{dp}$. Thus, the Lebesgue measure of this set is zero and almost all $\theta\in\bbr\times\bbc^{dp}$ satisfy the condition on the coefficients $d^*_i(\la_i)$ in Theorem~\ref{identCAR}. For instance, in the setting of Example~\ref{varioaxes}, $d_2^*(\la_{21})=0$ if and only if $\left(\la_{11}^2+3 \la_{11} \la_{12}+\la_{12}^2-\la_{21}^2\right)=0$.
\item The condition $-\pi/\Del\leq\Im(\la)<\pi/\Del$ is necessary due to the complex periodicity of the exponential function. In time series analysis this problem is associated with the \emph{aliasing effect}, i.e., the emergence of redundancies when sampling the process (cf. e.g. Section~3.4 and Assumption~C5 in \citet{Schlemm12}).\halmos
\end{enumerate}
\end{Remark}

Having established identifiability for CAR$(p)$ random fields, we now turn to CARMA$(p,q)$ random fields. For classical CARMA processes on $\bbr$ it is commonly known that one needs to impose at least conditions like $b_0\geq0$ and \emph{invertibility} in order to identify CARMA parameters from the second-order structure (i.e. either autocovariance, spectral density or variogram). For instance, if we consider the spectral density
\begin{equation*}
f(\om)=\frac{1}{2\pi}\frac{ |b(\bi\om)|^2}{ |a(\bi\om)|^2},\quad \om\in\bbr,
\end{equation*}
of a CARMA$(p,q)$ process with AR polynomial $a(\cdot)$ and MA polynomial $b(\cdot)$, then the numerator of $f$ yields the polynomial $b(z)b(-z)$. For every root $\la$ of $b(z)b(-z)$, $-\la$ is also a root, making it impossible to recover $b(\cdot)$ from $n(\cdot)$. Therefore, assuming invertibility, i.e., the condition that every root of $b(\cdot)$ has a negative real part, is necessary to determine the MA polynomial $b(\cdot)$ uniquely. However, this reasoning cannot be carried over to the causal CARMA random field since two additional obstacles occur: first, the spectral density $f$ of Proposition~\ref{specden} is now a multi-parameter function and, second, it is in general not separable, i.e., it cannot be written as a product of the form $f(\om)=f_1(\om_1)\cdots f_d(\om_d)$. Therefore, we cannot iterate the previous argument to each dimension. Also, the roots of the numerator $Q(z)Q(-z)$ are not discrete points in $\bbc$ anymore but, more generally, form algebraic varieties in $\bbc^d$. This makes it harder to formulate a similar condition as invertibility for the multi-parameter case. However, as we shall see by the end of this section, an invertibility condition is in fact not necessary. 
In order to show identifiability of CARMA random fields, we study the algebraic properties of the variogram and start with the following result.

\begin{Theorem}
Suppose that $(Y(t))_{t\in\bbr^d}$ is a CARMA$(p,q)$ random field such that Assumption~\ref{assump} holds true, $\kappa_2=1$, all eigenvalues $\la$ of $A_1,...,A_d$ satisfy $-\pi/\Del\leq\Im(\la)<\pi/\Del$ and all coefficients $d^*_i(\la_i)$ in Lemma~\ref{lem1} are nonzero. Further assume that the set $S=\{ \psi(j\Del\be_i)\colon i=1,...,d;j=0,...,2p+1 \}$ of variogram ordinates is given. Then there are at most $2^p$ different parameter values for $\theta$ which generate $S$.
\end{Theorem}

\begin{Proof}
Analogously to Theorem~\ref{identCAR}, we can determine all eigenvalues $\la$ of $A_1,...,A_d$ from the set $S$.
It remains to show that only finitely many vectors $b$ can generate $S$. By Lemma~\ref{lem1} and Assumption~\ref{assump}, we can solve Equations~\eqref{help12} for all coefficients $d^*_i(\la_i)$. By Theorem~4.1. in \cite{Pham18b} we have that $b$ has to satisfy the equations
\begin{equation}\label{eq1}
d^*_i(\la_i)=b^\top M(i,\la_i)b,
\end{equation}
where $i=1,...,d$, $\la_i$ is an eigenvalue of $A_i$ and $M(i,\la_i)$ are matrices that only depend on the (known) eigenvalues of $A_1,...,A_d$. That is, we are given $pd$ quadratic equations in $q+1$ unknowns $b_0,...,b_{q}$. Assumption~\ref{assump} and Bézout's theorem (see e.g. Theorem~18.3 in \cite{Harris92}) conclude the proof.

\halmos
\end{Proof}

The previous theorem shows that every fiber of the mapping $\Theta\ni\theta\mapsto S$ is finite. This property is also called \emph{algebraic identifiability} (cf. Section~1 in \cite{Amendola18}). 
To obtain statistical identifiability, we explicitly compute the variogram coefficients $d^*_i(\la_i)$ in Equation~\eqref{eq1}, which yields a polynomial system in terms of the CARMA parameters. One has then to show that this system has a unique solution. We demonstrate our method for the CARMA$(2,1)$ case and show how it can be applied to higher $(p,q)$.
\begin{Proposition}\label{identCARMA(2,1)}
Let $(Y(t))_{t\in\bbr^2}$ be a CARMA$(2,1)$ random field such that Assumption~\ref{assump} holds true, $\kappa_2=1$, $b_0\geq 0$, all eigenvalues $\la$ of $A_1,A_2$ satisfy $-\pi/\Del\leq\Im(\la)<\pi/\Del$ and all coefficients $d^*_i(\la_i)$ in Lemma~\ref{lem1} are nonzero. Furthermore, assume the additional condition $\la_{11}\la_{12}\neq\la_{21}\la_{22}$. Then $\theta$ is uniquely determined by $\{ \psi(j\Del\be_i)\colon i=1,2;j=0,...,5 \}$.
\halmos
\end{Proposition}
\begin{Proof}
First of all, the eigenvalues $\la_{11},\la_{12},\la_{21},\la_{22}$ and coefficients $d_1^*(\la_{11}),d_1^*(\la_{12}),d_2^*(\la_{21})$ and $d_2^*(\la_{22})$ can be recovered exactly as in Theorem~\ref{identCAR}. Hence, we only have to determine the parameters $b_0$ and $b_1$. Using the formulae of Example~\ref{varioaxes}, it is an easy task to verify the equations
\begin{align*}
d_1^*(\la_{11})&=\frac{ (b_0+b_1 \la_{11}) \left(b_0 \left(2 \la_{11} \la_{12}+\la_{12}^2+\la_{21}
\la_{22}\right)+b_1 \la_{11} \left(\la_{12}^2-\la_{21} \la_{22}\right)\right)}{4 \la_{11} \la_{21} \la_{22} (\la_{12}-\la_{11})
(\la_{11}+\la_{12}) (\la_{21}+\la_{22})},\\
d_1^*(\la_{12})&=\frac{ (b_0+b_1 \la_{12}) \left(b_0 \left(\la_{11}^2+2 \la_{11} \la_{12}+\la_{21} \la_{22}\right)+b_1 \la_{12}
\left(\la_{11}^2-\la_{21} \la_{22}\right)\right)}{4 \la_{12} \la_{21} \la_{22} (\la_{11}-\la_{12}) (\la_{11}+\la_{12})
(\la_{21}+\la_{22})},\\
d_2^*(\la_{21})&=\frac{ b_0^2 \left(\la_{11}^2+3 \la_{11} \la_{12}+\la_{12}^2-\la_{21}^2\right)+2 b_0 b_1 \la_{11} \la_{12}
(\la_{11}+\la_{12})+b_1^2 \la_{11} \la_{12} \left(\la_{11} \la_{12}-\la_{21}^2\right)}{4 \la_{11} \la_{12} \la_{21}
(\la_{11}+\la_{12}) (\la_{22}-\la_{21}) (\la_{21}+\la_{22})},\\
d_2^*(\la_{22})&=\frac{ b_0^2 \left(\la_{11}^2+3 \la_{11}
\la_{12}+\la_{12}^2-\la_{22}^2\right)+2 b_0 b_1 \la_{11} \la_{12} (\la_{11}+\la_{12})+b_1^2 \la_{11} \la_{12} \left(\la_{11}
\la_{12}-\la_{22}^2\right)}{4 \la_{11} \la_{12} \la_{22} (\la_{11}+\la_{12}) (\la_{21}-\la_{22}) (\la_{21}+\la_{22})}.
\end{align*}
We have to show that $b_0$ and $b_1$ are identifiable from this system, where $\la_{11},\la_{12},\la_{21},\la_{22}$ and $d_1^*(\la_{11}),d_1^*(\la_{12}),d_2^*(\la_{21}),d_2^*(\la_{22})$ are known. It therefore suffices to consider all four numerators and show that the system
\begin{align}
\notag& (\bar b_0+\bar b_1 \la_{11}) \left(\bar b_0 \left(2 \la_{11} \la_{12}+\la_{12}^2+\la_{21}
\la_{22}\right)+\bar b_1 \la_{11} \left(\la_{12}^2-\la_{21} \la_{22}\right)\right)\\
\notag&\quad= (b_0+b_1 \la_{11}) \left(b_0 \left(2 \la_{11} \la_{12}+\la_{12}^2+\la_{21}
\la_{22}\right)+b_1 \la_{11} \left(\la_{12}^2-\la_{21} \la_{22}\right)\right),\\
\notag& (\bar b_0+\bar b_1 \la_{12}) \left(\bar b_0 \left(\la_{11}^2+2 \la_{11} \la_{12}+\la_{21} \la_{22}\right)+\bar b_1 \la_{12}
\left(\la_{11}^2-\la_{21} \la_{22}\right)\right)\\
\notag&\quad= (b_0+b_1 \la_{12}) \left(b_0 \left(\la_{11}^2+2 \la_{11} \la_{12}+\la_{21} \la_{22}\right)+b_1 \la_{12}
\left(\la_{11}^2-\la_{21} \la_{22}\right)\right),\\
\notag& \bar b_0^2 \left(\la_{11}^2+3 \la_{11} \la_{12}+\la_{12}^2-\la_{21}^2\right)+2 \bar b_0 \bar b_1 \la_{11} \la_{12}
(\la_{11}+\la_{12})+\bar b_1^2 \la_{11} \la_{12} \left(\la_{11} \la_{12}-\la_{21}^2\right)\\
\notag&\quad= b_0^2 \left(\la_{11}^2+3 \la_{11} \la_{12}+\la_{12}^2-\la_{21}^2\right)+2 b_0 b_1 \la_{11} \la_{12}
(\la_{11}+\la_{12})+b_1^2 \la_{11} \la_{12} \left(\la_{11} \la_{12}-\la_{21}^2\right),\\
\notag& \bar b_0^2 \left(\la_{11}^2+3 \la_{11}
\la_{12}+\la_{12}^2-\la_{22}^2\right)+2 \bar b_0 \bar b_1 \la_{11} \la_{12} (\la_{11}+\la_{12})+\bar b_1^2 \la_{11} \la_{12} \left(\la_{11}
\la_{12}-\la_{22}^2\right)\\
\label{help14}&\quad= b_0^2 \left(\la_{11}^2+3 \la_{11}
\la_{12}+\la_{12}^2-\la_{22}^2\right)+2 b_0 b_1 \la_{11} \la_{12} (\la_{11}+\la_{12})+b_1^2 \la_{11} \la_{12} \left(\la_{11}
\la_{12}-\la_{22}^2\right),
\end{align}
implies $\bar b_0= b_0$ and $\bar b_1= b_1$, where we assume that $\bar b_0$ is non-negative and $\bar b_1\neq0$. Defining the variables
\begin{equation}\label{eq2}
x_1=\bar b_0^2- b_0^2,\quad x_2=\bar b_0 \bar b_1 - b_0 b_1,\quad x_3=\bar b_1^2- b_1^2,
\end{equation}
we find the equivalent linear system
\begin{equation}\label{linsystem}
\begin{pmatrix}
2\la_{11}\la_{12}+\la_{12}^2+\la_{21}\la_{22}	&2\la_{11}^2\la_{12}+2\la_{11}\la_{12}^2	&\la_{11}^2\la_{12}^2-\la_{11}^2\la_{21}\la_{22}\\
2\la_{11}\la_{12}+\la_{11}^2+\la_{21}\la_{22}	&2\la_{11}^2\la_{12}+2\la_{11}\la_{12}^2	&\la_{11}^2\la_{12}^2-\la_{12}^2\la_{21}\la_{22}\\
\la_{11}^2+3\la_{11}\la_{12}+\la_{12}^2-\la_{21}^2	&2\la_{11}^2\la_{12}+2\la_{11}\la_{12}^2	&\la_{11}^2\la_{12}^2-\la_{11}\la_{12}\la_{21}^2\\
\la_{11}^2+3\la_{11}\la_{12}+\la_{12}^2-\la_{22}^2	&2\la_{11}^2\la_{12}+2\la_{11}\la_{12}^2	&\la_{11}^2\la_{12}^2-\la_{11}\la_{12}\la_{22}^2
\end{pmatrix}x=0,
\end{equation}
with $x^\top=(x_1,x_2,x_3)^\top$. This system has the unique solution $x=0$ if and only if at least one of the four $3\times3$-minors
\begin{align*}
2 \la_{11} \la_{12} \la_{21} (\la_{11}-\la_{12}) (\la_{11}+\la_{12})^2 (\la_{21}+\la_{22}) (\la_{11} \la_{12}-\la_{21} \la_{22}),\\
2 \la_{11} \la_{12} \la_{22} (\la_{11}-\la_{12}) (\la_{11}+\la_{12})^2 (\la_{21}+\la_{22}) (\la_{11} \la_{12}-\la_{21} \la_{22}),\\
-2 \la_{11}^2 \la_{12} (\la_{11}+\la_{12})^2 (\la_{21}-\la_{22}) (\la_{21}+\la_{22}) (\la_{11} \la_{12}-\la_{21} \la_{22}),\\
-2 \la_{11} \la_{12}^2 (\la_{11}+\la_{12})^2 (\la_{21}-\la_{22}) (\la_{21}+\la_{22}) (\la_{11} \la_{12}-\la_{21} \la_{22}),
\end{align*}
is not zero. However, this is equivalent to the condition $\la_{11}\la_{12}\neq\la_{21}\la_{22}$.
Hence, by our assumptions we can indeed conclude that $x=0$, which yields $\bar b_0= b_0$ and $\bar b_1= b_1$.
\halmos
\end{Proof}
\begin{Remark}
\begin{enumerate}
\item Note that in Proposition~\ref{identCARMA(2,1)} we have not used the full variogram but only values on the principal axes. Working with the full variogram, we are able to dispose of the condition $\la_{11}\la_{12}\neq\la_{21}\la_{22}$. However, imposing this weak condition has the advantage that we do not have to estimate the full variogram and the set of parameters which satisfy $\la_{11}\la_{12}=\la_{21}\la_{22}$ is a Lebesgue null set in $\bbc^4$.
\item In the setting of Proposition~\ref{identCARMA(2,1)} the condition $\la_{11}\la_{12}\neq\la_{21}\la_{22}$ is not only sufficient but also necessary. For instance, if we choose $\la_{11}=\la_{21}=-2$ and $\la_{12}=\la_{22}=-6$, then both pairs $(b_0,b_1)=(2,4)$ and $(b_0,b_1)=(20/\sqrt{7},9/\sqrt{7})$ will generate the same variogram on the principal axes. Hence, in this case we do not have identifiability of the model parameters.

\halmos
\end{enumerate}
\end{Remark}

In a similar fashion we can show the following result. Since all factors in \eqref{help16} are nonzero, there is no extra condition like $\la_{11}\la_{12}\neq\la_{21}\la_{22}$ needed as in Proposition~\ref{identCARMA(2,1)}.

\begin{Proposition}\label{identCARMA(3,1)}
Let $(Y(t))_{t\in\bbr^2}$ be a CARMA$(3,1)$ random field such that Assumption~\ref{assump} holds true, $\kappa_2=1$, $b_0\geq 0$, all eigenvalues $\la$ of $A_1,A_2$ satisfy $-\pi/\Del\leq\Im(\la)<\pi/\Del$ and all coefficients $d^*_i(\la_i)$ in Lemma~\ref{lem1} are nonzero. Then $\theta$ is uniquely determined by $\{ \psi(j\Del\be_i)\colon i=1,2;j=0,...,7 \}$.
\halmos
\end{Proposition}

\begin{Proof}
The assertion can be proven analogously to the proof of Proposition~\ref{identCARMA(2,1)}. We therefore only highlight the difference. Instead of 4 we have 6 different $d^*_i(\la_i)$ in this case. Defining $x_1,x_2,x_3$ as before, we obtain a linear system of size $6\times3$ similar to Equation~\eqref{linsystem}. The system matrix has $\binom{6}{3}=20$ different $3\times3$-minors, one of which is 
\begin{align}\notag
&-6 \la_{11}^2 \la_{12}^2 \la_{13}^2 (\la_{11}+\la_{12})^2 (\la_{11}+\la_{13})^2 (\la_{12}+\la_{13})^2 (\la_{21}-\la_{22}) (\la_{21}+\la_{22})\\
\label{help16}&\quad\times(\la_{21}-\la_{23}) (\la_{21}+\la_{23}) (\la_{22}-\la_{23}) (\la_{22}+\la_{23}).
\end{align}
This minor is always nonzero under our assumptions. Thus, we conclude $\bar b_0= b_0$ and $\bar b_1= b_1$.

\halmos
\end{Proof}

The method used to show identifiability for CARMA$(2,1)$ and CARMA$(3,1)$ random fields on $\bbr^2$ relied on the definition of appropriate variables $x_1,x_2,x_3$ in Equation~\eqref{eq2} and a system of $pd=4$ equations in the first case and $pd=6$ equations in the second case. Both systems have a unique solution provided that at least one of the minors of the coefficient matrix is nonzero. In the first case $4$ minors of the $4\times3$ coefficient matrix had to be considered and in the second case $20$ of the $6\times3$ coefficient matrix. The complexity of this method becomes too high to consider higher order models. Moreover, we have observed that for the CARMA$(3,2)$ model on $\bbr^2$ the method fails, since the determinant of the corresponding $6\times6$ coefficient matrix is always zero. However, this does not prevent  parameter identifiability, since -- as we note from Equation~\eqref{eq2} -- the components of the vector $x$ display algebraic dependencies, that is, the variables of the corresponding linear systems are not independent.

As an alternative to the substitution \eqref{eq2}, we can find a solution to the original system of $pd$ quadratic equations \eqref{eq1} for the $q+1$ variables $b_0,\dots,b_q$ directly taking resort to representations via Gröbner bases (see e.g. Chapter~2 of \citet{Cox15}). As a test case we have replicated Proposition~\ref{identCARMA(2,1)} using the software \texttt{Mathematica}, where the $pd=4$ quadratic equations in \eqref{help14} were transformed to an equivalent system of $48$ polynomial equations. From these we could read off $\bar b_0= b_0$ and $\bar b_1= b_1$ immediately and again obtain identifiability. 

Note that in Propositions~\ref{identCARMA(2,1)} and \ref{identCARMA(3,1)} we have not assumed any extra conditions on $b$ except for $b_0\geq0$. In particular, it is not necessary to impose an analogous condition to invertibility in order to achieve identifiability. This illustrates a fundamental difference between CARMA processes and CARMA random fields with $d\geq2$.

\section{Simulation of CARMA random fields on a lattice}\label{sectionsimu}

In this section we develop two numerical simulation schemes for the causal CARMA random field. One is designed for compound Poisson noise and the other one for general Lévy noise. In both cases, we simulate on a lattice $L=\{ \Del,...,N\Del \}^d$ with fixed $\Del>0$ and $N\in\bbn$. Techniques for simulating on more general lattices are discussed as well.

\subsection{Compound Poisson noise}
The homogeneous Lévy basis $\La$ is assumed to be compound Poisson in this subsection. That is, the characteristic triplet of $\La$ satisfies $\beta=c\int_{(-1,1)}x\, F(\dd x)$, $\sigma=0$ and $\nu=cF$, where $c>0$ is the intensity parameter and $F$ is a probability measure on $\bbr$ (cf. Section~1.2.4 in \citet{Applebaum09}). As a consequence, the resulting CARMA random field $(Y(t))_{t\in\bbr^d}$ in Equation~\eqref{convo} can be represented as
\begin{equation*}
Y(t)=\sum_{j\in\bbn}g(t-s_j)W(s_j),\quad t\in\bbr^d,
\end{equation*}
where $s_j$ are the locations of the countably many Lévy jumps of $\La$ and the i.i.d. $W(s_j)$ are the heights of the Lévy jumps, distributed according to $F$. Restricted on a compact domain $D\subset\bbr^d$, there are only finitely many jumps of $\La$ and their number $N_J$ follows a Poisson distribution with intensity $c\Leb(D)$. Conditionally on the value of $N_J$, the jump positions are independently and uniformly distributed on $D$. This motivates us to approximate $Y$ with
\begin{equation*}
Y_{S1}(t)=\sum_{j=1}^{N_J}g(t-s_j)W(s_j),\quad t\in\bbr^d,\quad s_j\in D.
\end{equation*}
The random field $Y_{S1}$ has the alternative representation
\begin{equation}
\label{convoS1} Y_{S1}(t)=(g*\La_{S1})(t):=\int_{\bbr^d} g(t-s) \,\La_{S1}(\dd s),\quad t\in\bbr^d,
\end{equation}
with $\La_{S1}(\dd s)=\bone_D(s)\La(\dd s)$, hence it arises by truncating the Lévy basis $\La$. The advantage of $Y_{S1}$ is that we can simulate it exactly.
For the simulation algorithm we choose $D=[-M,M]^d$ with a sufficiently large $M>0$ such that $L\subset D$.

\pagebreak

\begin{Algorithm}~\label{algo1}
\begin{enumerate}
\item Input: $g,F,c,M,N,\Del$ such that $\{ \Del,...,N\Del \}^d=L\subset D=[-M,M]^d$
\item Draw $N_J$ from a Poisson distribution with intensity $c\Leb(D)=c(2M)^d$.
\item Draw $s_1,...,s_{N_J}$ independently and uniformly distributed on $D=[-M,M]^d$.
\item For each $s_j$, $j=1,...,N_J$, draw $W(s_j)$ independently from the distribution $F$.
\item For each $t\in L$, compute $Y_{S1}(t)=\sum_{j=1}^{N_J}g(t-s_j)W(s_j)$. \label{help4}
\item Output: $Y_{S1}(t),\quad t\in L=\{ \Del,...,N\Del \}^d$
\end{enumerate}
\halmos
\end{Algorithm}
In order to assess the accuracy of this approximation algorithm, we determine its mean squared error. Note that as the simulation of $Y_{S1}$ is exact, we only have to consider the approximation error between $Y$ and $Y_{S1}$. Moreover, we show that the simulated random field $Y_{S1}(t)$ converges for fixed $t\in L$ both in $L^2$ and almost surely to the underlying true random field $Y(t)$ as the truncation parameter $M$ tends to infinity.
\begin{Theorem}\label{TSim1}
Suppose that $(Y(t))_{t\in\bbr^d}$ is a CARMA$(p,q)$ random field such that Assumption~\ref{assump} holds true and $\La$ is compound Poisson with characteristic triplet $(c\int_{(-1,1)}x\, F(\dd x),0,cF)$, where $c>0$ and $F$ is a probability distribution. Then the mean squared error of Algorithm~\ref{algo1} satisfies
\begin{align}
\notag\max_{t\in L}\bbe\Big[\big(Y(t)-Y_{S1}(t)\big)^2\Big]
&= \kappa_2\sum_{\la_1}\cdots\sum_{\la_d}\sum_{\la_1'}\cdots\sum_{\la_d'} d(\la_1,...,\la_d)d(\la_1',...,\la_d')\\
\notag&\quad\times\Bigg[\frac{1}{|\la_1|+|\la_1'|}\cdots\frac{1}{|\la_d|+|\la_d'|}
-\frac{1-\ee^{(\la_1+\la_1')M}}{|\la_1|+|\la_1'|}\cdots\frac{1-\ee^{(\la_d+\la_d')M}}{|\la_d|+|\la_d'|}
\Bigg]\\
\label{help1}&=\calo(\ee^{-2|\la_{max}|M}), \quad M\to\infty,
\end{align}
where the coefficients $d(\cdot)$ are the same as in Equation~\eqref{kernel2}, both $\sum_{\la_i}$ and $\sum_{\la_i'}$ denote the sum over distinct eigenvalues of $A_i$ for $i=1,...,d$ and $\la_{max}:=\max\{ \Re(\la)\colon \la \text{ is eigenvalue of } A_i,i=1,...,d \}$.

Furthermore, $Y_{S1}(t)$ converges to $Y(t)$ in $L^2$ and almost surely as $M\to\infty$ for every $t\in L=\{ \Del,...,N\Del \}^d$.
\end{Theorem}

\begin{Proof}
By the properties of Lévy bases and Equations~\eqref{convo} and \eqref{convoS1} we observe that
\begin{align}
\notag\max_{t\in L}\bbe\Big[\big(Y(t)-Y_{S1}(t)\big)^2\Big] &= \max_{t\in L}\bbe\Big[\big( \int_{[-\infty\be,t]/[-M\be,t]} g(t-s) \,\La(\dd s) \big)^2\Big]\\
\notag&=\max_{t\in L} \kappa_2\left(\int_{\bbr_+^d} g^2(s)\,\dd s-\int_{[0,t+M\be]} g^2(s)\,\dd s\right)\\
\label{help2}&=\kappa_2\int_{\bbr_+^d} g^2(s)\,\dd s - \kappa_2\int_{[0,M]^d} g^2(s)\,\dd s = \bbe\Big[\big(Y(0)-Y_{S1}(0)\big)^2\Big],
\end{align}
where in the first equation we have taken into account that the kernel $g$ contains the indicator function $\bone_{\{s\geq0\}}$. In addition, Equation~\eqref{kernel2} implies
\begin{equation*}
g^2(s)=\sum_{\la_1}\cdots\sum_{\la_d}\sum_{\la_1'}\cdots\sum_{\la_d'} d(\la_1,...,\la_d)d(\la_1',...,\la_d') \ee^{(\la_1+\la_1') s_1}\cdots \ee^{(\la_d+\la_d') s_d}\bone_{\{s\geq0\}},\quad s\in\bbr^d.
\end{equation*}
Plugging this into \eqref{help2}, we arrive at Equation~\eqref{help1}, which in turn shows that $Y_{S1}(t)$ converges to $Y(t)$ in $L^2$ for every $t\in L=\{ \Del,...,N\Del \}^d$. It remains to show that the convergence also holds almost surely. Owing to Chebyshev's inequality we have for each $t\in L=\{ \Del,...,N\Del \}^d$ that
\begin{equation*}
\sum_{M=1}^\infty \bbp\left[|Y(t)-Y_{S1,M}(t)|\geq \frac{1}{M}\right]\leq \sum_{M=1}^\infty M^2\bbe[|Y(t)-Y_{S1,M}(t)|^2],
\end{equation*}
where we explicitly include the input parameter $M$ into the subscript of $Y_{S1,M}(t)$. The right-hand side of the latter inequality is finite due to Equation~\eqref{help1}. Finally, the assertion follows from  the Borel-Cantelli lemma.
\halmos
\end{Proof}

\begin{Remark}
\begin{enumerate}
\item Algorithm~\ref{algo1} can also be applied to pure-jump Lévy bases if small jumps are truncated. This technique has been analyzed in detail in Section~3 of \citet{Chen16} for the simulation of stochastic Volterra equations in space--time. Furthermore, Section~4 of \cite{Chen16} considers a simulation technique which is based on series representations for Lévy bases (see also \citet{Rosinski01}). However, we do not pursue this direction. Instead, in the next subsection we consider a method which are not restricted to pure-jump Lévy bases, easy to implement and sufficient for our simulation study in Section~\ref{simstudy}.
\item One can readily replace $L=\{ \Del,...,N\Del \}^d$ in step~(\ref{help4}) of Algorithm~\ref{algo1} with any finite subset of points in $\bbr^d$. Algorithm~\ref{algo1} is not restricted to simulation on lattices.
\halmos
\end{enumerate}
\end{Remark}

\subsection{General Lévy noise}

Algorithm~\ref{algo1} is not suitable for CARMA random fields driven by general Lévy bases since a drift or a Gaussian part may be part of the noise. A different way to approximate a CARMA random field $(Y(t))_{t\in\bbr^d}$ is to discretize and truncate the stochastic integral in Equation~\eqref{convo}. Introducing a truncation parameter $M\in\bbn$, we first replace the integral in \eqref{convo} by
\begin{equation*}
\int_{[t-\Del M\be,t]} g(t-s) \,\La(\dd s),\quad t\in\bbr^d.
\end{equation*}
By discretization of this integral we obtain the sum 
\begin{align}
\notag Y_{S2}(t):&=\sum_{s\in [t-\Del M\be,t]\cap\Del\bbz^d} g(t-s)Z(s)\\
\label{help6}&=\sum_{s\in \{ 0,\Del,...,M\Del \}^d} g(s)Z(t-s),\quad t\in\bbr^d.
\end{align}
Here, the random field $Z$ represents spatial increments of $\La$, or more precisely
\begin{equation*}
Z(t):=\La([t-\Del,t]),\quad t\in\bbr^d.
\end{equation*}
This approach has also been applied in \cite{Nguyen15} to simulate the so-called $OU_\wedge$ process. Since we evaluate $Y_{S2}$ only on the lattice $L=\{ \Del,...,N\Del \}^d$, we actually simulate a discrete-parameter moving average random field of finite order driven by i.i.d. spatial noise as given in \eqref{help6}. The set $\{g(s)\colon s\in \{ 0,\Del,...,M\Del \}^d\}$ plays the role of the moving average coefficients and $Y_{S2}$ can be simulated exactly. Furthermore, it is easy to check that the random field $Y_{S2}$ also has the representation
\begin{equation*}
\label{convoS2} Y_{S2}(t)=(g_{S2}*\La)(t):=\int_{\bbr^d} g_{S2}(t-s) \,\La(\dd s),\quad t\in\bbr^d,
\end{equation*}
where the step function $g_{S2}$ is given by
\begin{equation}\label{stepfunc}
g_{S2}(s)=\sum_{j\in\{ 0,,...,M \}^d}b^\top\ee^{A_1 j_1\Del}\cdots\ee^{A_d j_d\Del}\be_p\bone_{[j\Del,(j+\be)\Del]}(s),\quad s\in\bbr^d.
\end{equation}
This allows us to observe that truncation and discretization of the stochastic integral in Equation~\eqref{convo} is in fact equivalent to truncation and discretization of the kernel $g$, which will be useful for establishing error bounds. We sum up the simulation scheme in the following algorithm, where $(\beta,\si^2,\nu)$ denotes the characteristic triplet of $\La$.
\begin{Algorithm}\label{algo2}~
\begin{enumerate}
\item Input: $g,(\beta,\si^2,\nu),M,N,\Del$
\item Compute $g(s)$ for $s\in \{ 0,\Del,...,M\Del \}^d$.
\item Draw $Z(s)$, $s\in\{ (1-M)\Del,...,N\Del \}^d$, independently from the infinitely divisible distribution with characteristics $(\Del^d\beta,\Del^d\si^2,\Del^d\nu)$.
\item For each $t\in L=\{ \Del,...,N\Del \}^d$, compute $Y_{S2}(t)=\sum_{s\in \{ 0,\Del,...,M\Del \}^d} g(s)Z(t-s)$.
\item Output: $Y_{S2}(t),\quad t\in L=\{ \Del,...,N\Del \}^d$
\end{enumerate}
\halmos
\end{Algorithm}
If we collect the $g(s)$ values from the second step of Algorithm~\ref{algo2} in an array $A_g$, and the $Z(s)$ values from the third step in an array $A_Z$, then the $Y_{S2}(s)$ values from the fourth step can be computed as the discrete convolution of the two arrays $A_g$ and $A_Z$. This can be carried out efficiently using the fast Fourier transform (FFT). In-built convolution commands using the FFT exist in computer softwares such as \texttt{R} or \texttt{Matlab}.

By approximating the CARMA random field $Y$ by $Y_{S2}$ we create two sources of error, one originates from the kernel truncation, the other one from the kernel discretization. A more detailed analysis yields the following result.
\begin{Theorem}\label{TSim2}
Suppose that $(Y(t))_{t\in\bbr^d}$ is a CARMA$(p,q)$ random field such that Assumption~\ref{assump} holds true. Then $Y_{S2}(t)$ converges to $Y(t)$ in $L^2$ as simultaneously $\Del\to0$ and $\Del M\to\infty$ for every $t\in L=\{ \Del,...,N\Del \}^d$.

Further, let $(\Del_k)_{k\in\bbn}$ and $(M_k)_{k\in\bbn}$ be two sequences satisfying $\Del_k=\calo(k^{-1-\eps})$ for some $\eps>0$ and $\Del_k M_k\to\infty$ as $k\to\infty$. Then $Y_{S2}(t)$ also converges to $Y(t)$ almost surely as $k\to\infty$ for every $t\in L=\{ \Del,...,N\Del \}^d$.
\end{Theorem}
\begin{Proof}
For notational convenience we assume that all eigenvalues of $A_1,...,A_d$ are real. The complex case can be shown analogously by similar arguments taking care of imaginary parts. The mean squared error is by stationarity for each $t\in L=\{ \Del,...,N\Del \}^d$ the same, namely
\begin{align}
\label{help17}&\bbe\Big[\big(Y(t)-Y_{S2}(t)\big)^2\Big] = \kappa_2\int_{\bbr^d} (g(s)-g_{S2}(s))^2\,\dd s\\
\notag&\quad\leq \kappa_2p^d\sum_{\la_1}\cdots\sum_{\la_d} d(\la_1,...,\la_d)^2\\
\notag&\quad\quad\times\int_{\bbr_+^d}\bigg(\ee^{\la_1 s_1}\cdots \ee^{\la_d s_d} 
- \sum_{j\in\{ 0,,...,M \}^d}\ee^{\la_1 j_1\Del}\cdots \ee^{\la_d j_d\Del}\bone_{[j\Del,(j+\be)\Del]}(s) \bigg)^2\,\dd s.
\end{align}
Here, we have used the inequality $(\sum_{j=1}^K a_j)^2\leq K\sum_{j=1}^K a_j^2$ with $a_1,...,a_K\in\bbr$. In order to evaluate the latter integral, we consider for fixed $\la_1,...,\la_d$ the identities
\begin{align*}
&\int_{\bbr_+^d}\ee^{2\la_1 s_1}\cdots \ee^{2\la_d s_d}\,\dd s = \frac{1}{-2\la_1}\cdots\frac{1}{-2\la_d},\\
&-2\int_{\bbr_+^d}\ee^{\la_1 s_1}\cdots \ee^{\la_d s_d}\sum_{j\in\{ 0,,...,M \}^d}\ee^{\la_1 j_1\Del}\cdots \ee^{\la_d j_d\Del}\bone_{[j\Del,(j+\be)\Del]}(s)\,\dd s\\
&\quad= -2\sum_{j\in\{ 0,,...,M \}^d}\frac{\ee^{\la_1(2j_1+1)\Del}-\ee^{2\la_1j_1\Del}}{\la_1}\cdots\frac{\ee^{\la_d(2j_d+1)\Del}-\ee^{2\la_dj_d\Del}}{\la_d},
\end{align*}
and
\begin{equation*}
\int_{\bbr_+^d}\sum_{j\in\{ 0,,...,M \}^d}\ee^{2\la_1 j_1\Del}\cdots \ee^{2\la_d j_d\Del}\bone_{[j\Del,(j+\be)\Del]}(s)\,\dd s = \sum_{j\in\{ 0,,...,M \}^d}\ee^{2\la_1 j_1\Del}\cdots \ee^{2\la_d j_d\Del}\Del^d.
\end{equation*}
Summing these up, we obtain
\begin{align*}
&\int_{\bbr_+^d}\bigg(\ee^{\la_1 s_1}\cdots \ee^{\la_d s_d} 
- \sum_{j\in\{ 0,,...,M \}^d}\ee^{\la_1 j_1\Del}\cdots \ee^{\la_d j_d\Del}\bone_{[j\Del,(j+\be)\Del]}(s) \bigg)^2\,\dd s\\
&\quad=\frac{1}{-2\la_1}\cdots\frac{1}{-2\la_d} + \sum_{j\in\{ 0,,...,M \}^d}\ee^{2\la_1 j_1\Del}\cdots \ee^{2\la_d j_d\Del}\left( \Del^d-2\frac{\ee^{\la_1\Del}-1}{\la_1}\cdots\frac{\ee^{\la_d\Del}-1}{\la_d} \right)\\
&\quad=f_1(\Del M,\Del),
\end{align*}
where the function $f_1$ is defined as
\begin{align*}
f_1(u,v)&:=\frac{1}{-2\la_1}\cdots\frac{1}{-2\la_d}\\
&\qquad+\frac{(\ee^{2\la_1(u+v)}-1)v}{\ee^{2\la_1v}-1}\cdots\frac{(\ee^{2\la_d(u+v)}-1)v}{\ee^{2\la_dv}-1} \left( 1-2\frac{\ee^{\la_1v}-1}{\la_1v}\cdots\frac{\ee^{\la_dv}-1}{\la_dv} \right).
\end{align*}
Additionally, we have the limits
\begin{align}
\label{help3}&\lim_{u\to\infty}f_1(u,v)=\frac{1}{-2\la_1}\cdots\frac{1}{-2\la_d} + \frac{-v}{\ee^{2\la_1v}-1}\cdots\frac{-v}{\ee^{2\la_dv}-1} \left( 1-2\frac{\ee^{\la_1v}-1}{\la_1v}\cdots\frac{\ee^{\la_dv}-1}{\la_dv} \right)\\
\notag&\quad=:f_2(v)=:\frac{1}{-2\la_1}\cdots\frac{1}{-2\la_d}+f_3(v),
\end{align}
and
\begin{equation}\label{help7}
\lim_{v\to 0}f_2(v)=0.
\end{equation}
Moreover, we observe that
\begin{equation*}
f_1(u,v)-f_2(v)=f_3(v)[(1-\ee^{2\la_1(u+v)})\cdots(1-\ee^{2\la_d(u+v)})-1].
\end{equation*}
For every $\eps>0$, the function $f_3(v)$ is bounded and continuous on $(0,\eps)$. We therefore arrive at
\begin{equation*}
\lim_{u\to\infty}\sup_{v\in(0,\eps)}|f_1(u,v)-f_2(v)|=0,
\end{equation*}
which shows that the convergence in Equation~\eqref{help3} is actually uniform in $v$. Combined with \eqref{help7}, this implies that
\begin{equation*}
\lim_{u\to\infty,v\to0}f_1(u,v)=0.
\end{equation*}
Hence, for every $t\in L=\{ \Del,...,N\Del \}^d$, $Y_{S2}(t)$ converges to $Y(t)$ in $L^2$ as simultaneously $\Del\to0$ and $\Del M\to\infty$.

As for the second part of our assertion, we note that if $\Del_k=\calo(k^{-1-\eps})$, then all $\la<0$ satisfy the inequality
\begin{equation}\label{help11}
\sum_{k=1}^\infty\left| \la-\frac{e^{\la\Del_k}-1}{\Del_k} \right|<\infty,
\end{equation}
which can be shown with the Taylor expansion of the exponential function. Defining
\begin{equation*}
A_{1,k}:=\frac{(\ee^{2\la_1(\Del_k M_k+\Del_k)}-1)\Del_k}{\ee^{2\la_1\Del_k}-1}\cdots\frac{(\ee^{2\la_d(\Del_k M_k+\Del_k)}-1)\Del_k}{\ee^{2\la_d\Del_k}-1} \left( 2-2\frac{\ee^{\la_1\Del_k}-1}{\la_1\Del_k}\cdots\frac{\ee^{\la_d\Del_k}-1}{\la_d\Del_k} \right)
\end{equation*}
and
\begin{equation*}
A_{2,k}:=\frac{1}{-2\la_1}\cdots\frac{1}{-2\la_d} - \frac{(\ee^{2\la_1(\Del_k M_k+\Del_k)}-1)\Del_k}{\ee^{2\la_1\Del_k}-1}\cdots\frac{(\ee^{2\la_d(\Del_k M_k+\Del_k)}-1)\Del_k}{\ee^{2\la_d\Del_k}-1}, 
\end{equation*}
Inequality~\eqref{help11} implies $\sum_{k=1}^\infty |A_{1,k}|<\infty$ and $\sum_{k=1}^\infty |A_{2,k}|<\infty$, and thus
\begin{equation*}
\sum_{k=1}^\infty |f_1(\Del_k M_k,\Del_k)|<\infty.
\end{equation*}
Finally, the almost sure convergence follows similarly as in the proof of Theorem~\ref{TSim1} by Chebyshev's inequality and the Borel-Cantelli lemma.
\halmos
\end{Proof}

\begin{Remark}
\begin{enumerate} 
\item Instead of simulating on the regular lattice $L=\{ \Del,...,N\Del \}^d$, one can easily adjust Algorithm~\ref{algo2} for simulating on the more general lattice $L=\{ \Del_1,...,N_1\Del_1 \}\times\cdots\times\{ \Del_d,...,N_d\Del_d \}$ with $\Del_1,...,\Del_d>0$ and $N_1,...,N_d\in\bbn$. 
\item In Section~4 of \cite{Pham18b} it was shown that under mild conditions every CARMA random field has a version which is càdlàg with respect to the partial order $\leq$. By inspection of Algorithm~\ref{algo1} and Algorithm~\ref{algo2} it is easy to see that both $Y_{S1}$ and $Y_{S2}$ are càdlàg as well.
\halmos
\end{enumerate}
\end{Remark}

\section{Simulation study}\label{simstudy}

We conduct a simulation study in order to assess the empirical quality of the WLS estimator of the previous section for finite samples. We use Algorithm~\ref{algo2} to simulate 500 paths of a CARMA$(2,1)$ random field on a two-dimensional grid. As CARMA parameters we take the estimates from Section~\ref{sectionCMB}, which are
\begin{equation*}
b_0=4.8940,\, b_1=-1.1432,\, \la_{11}=-1.7776,\, \la_{12}=-2.0948,\, \la_{21}=-1.3057,\, \la_{22}=-2.5142.
\end{equation*}
We take a Gaussian Lévy basis $\La$ with mean zero and variance one. In accordance with the parameter estimation in Section~\ref{sectionCMB}, we first choose $\Del=0.04$ for the grid size of Algorithm~\ref{algo2}, $M=400$ for the truncation parameter and $N^2=1000^2$ for the number of points for each path. However, this choice results in relatively high approximation errors, yielding only poor parameter estimates. By choosing a higher truncation parameter $M$ and a smaller grid size $\Del$, the step function $g_{S2}$ in \eqref{stepfunc} approximates the CARMA kernel $g$ in \eqref{kernel} better, which by \eqref{help17} also reduces the approximation error of the CARMA random field. We therefore decide to simulate on a finer grid with $\Del=0.01$, $M=600$ and $N^2=4000^2$. After simulation we save only every fourth point in each of the two axes directions of $\bbr^2$ in order to be back in the setting of Section~\ref{sectionCMB} with $\Del=0.04$ and $N^2=1000^2$ points per path.

Having simulated the CARMA random fields on a grid, we proceed by estimating the variogram using the variogram estimator of Section~\ref{sectionempvario}. We calculate the empirical variogram at $K=100$ different lags, namely
\begin{equation}\label{varioest}
\{ \psi_{N}^*(j\Del\be_i)\colon i=1,2;j=1,...,50 \}.
\end{equation}
These lags lie on the principal axes of $\bbr^2$ and are by Proposition~\ref{identCARMA(2,1)} sufficient to identify the CARMA$(2,1)$ parameters. In the final step we estimate the CARMA parameter vector $\theta$ with the WLS estimator $\theta_{N}^*$ given in \eqref{WLSest}. We consider the following choices for weights and number of lags used.

\vspace{\baselineskip}
\noindent \textbf{Case 1:}
\begin{equation}\label{case1}
\theta_{N}^*:=\argmin_{\theta\in\Theta}\Bigg\{ \sum_{\substack{j=1,...,50\\i=1,2}} w_j\left(\psi_{N}^*(j\Del\be_i)-\psi_{\theta}(j\Del\be_i)\right)^2 \Bigg\},\quad
w_j=\left(\frac{0.1(j-1)+50-j}{49}\right)^2.
\end{equation}
\textbf{Case 2:}
\begin{equation}\label{case2}
\theta_{N}^*:=\argmin_{\theta\in\Theta}\Bigg\{ \sum_{\substack{j=1,...,25\\i=1,2}} w_j\left(\psi_{N}^*(j\Del\be_i)-\psi_{\theta}(j\Del\be_i)\right)^2 \Bigg\},\quad
w_j=\left(\frac{0.1(j-1)+25-j}{24}\right)^2.
\end{equation}
\textbf{Case 3:}
\begin{equation}\label{case3}
\theta_{N}^*:=\argmin_{\theta\in\Theta}\Bigg\{ \sum_{\substack{j=1,...,50\\i=1,2}} w_j\left(\psi_{N}^*(j\Del\be_i)-\psi_{\theta}(j\Del\be_i)\right)^2 \Bigg\},\quad
w_j=\ee^{j\Del}.
\end{equation}
\textbf{Case 4:}
\begin{equation}\label{case4}
\theta_{N}^*:=\argmin_{\theta\in\Theta}\Bigg\{ \sum_{\substack{j=1,...,25\\i=1,2}} w_j\left(\psi_{N}^*(j\Del\be_i)-\psi_{\theta}(j\Del\be_i)\right)^2 \Bigg\},\quad
w_j=\ee^{j\Del}.
\end{equation}
Cases 1 and 2 apply quadratically decreasing weights while Cases 3 and 4 apply exponentially decreasing weights.
The compact parameter space $\Theta$ is chosen to be $\Theta=[0,10]\times[-10,10]\times[-10,0]^4$ which contains the true parameter vector $\theta_0=(b_0, b_1, \la_{11}, \la_{12}, \la_{21}, \la_{22})$. For minimization of the objective function we use the command \texttt{DEoptim} of the \texttt{R} package \textbf{DEoptim} which implements the differential evolution algorithm (for more details see \cite{Mullen11}). This algorithm has the advantage that we do not need an initial value for the optimization procedure. Instead, one can directly hand over the parameter space $\Theta$ as an input. The output of \texttt{DEoptim} itself is then used as an initial value for the standard \texttt{R} command \texttt{optim}. The summary of the estimation results are given in Tables~\ref{summary1} to \ref{summary4} below.

Recall that in our parametrization $b_0$ actually plays the role of the white noise standard deviation. Comparing Table~\ref{summary1}
with \ref{summary2} and Table~\ref{summary3}
with \ref{summary4}, we observe that using $K=50$ instead of $K=100$ lags generally reduces the standard deviation (Std) but increases the bias for most of the estimators. This indicates a typical variance-bias trade-off subject to the number of lags used. Moreover, we find that using exponential weights as in \eqref{case3} and \eqref{case4} increases the standard deviation and the root mean squared error (RMSE) for all components of $\theta_{N}^*$.

According to Theorem~\ref{thmWLS}, the asymptotic properties of the WLS estimator $\theta_{N}^*$ does not depend on the distribution of the Lévy basis $\La$. To examine this statement for finite samples, we repeat the procedure above with variance gamma noise. More precisely, we simulate 500 independent CARMA$(2,1)$ paths driven by a variance gamma basis $\La$ with mean zero and variance one, compute the empirical variogram as in \eqref{varioest} and estimate the CARMA parameters as in Cases~1 to 4. The results are summarized in Tables~\ref{summary5} to ~\ref{summary8}. Comparing the RMSEs in Tables~\ref{summary1} to \ref{summary8}, we observe that the WLS estimation is slightly but not significantly better for the variance gamma case than for the Gaussian case.

\section{Application to cosmic microwave background data}\label{sectionCMB}

\begin{figure}[ht]
  \begin{minipage}[b]{0.5\linewidth}
    \centering
    \includegraphics[width=\linewidth]{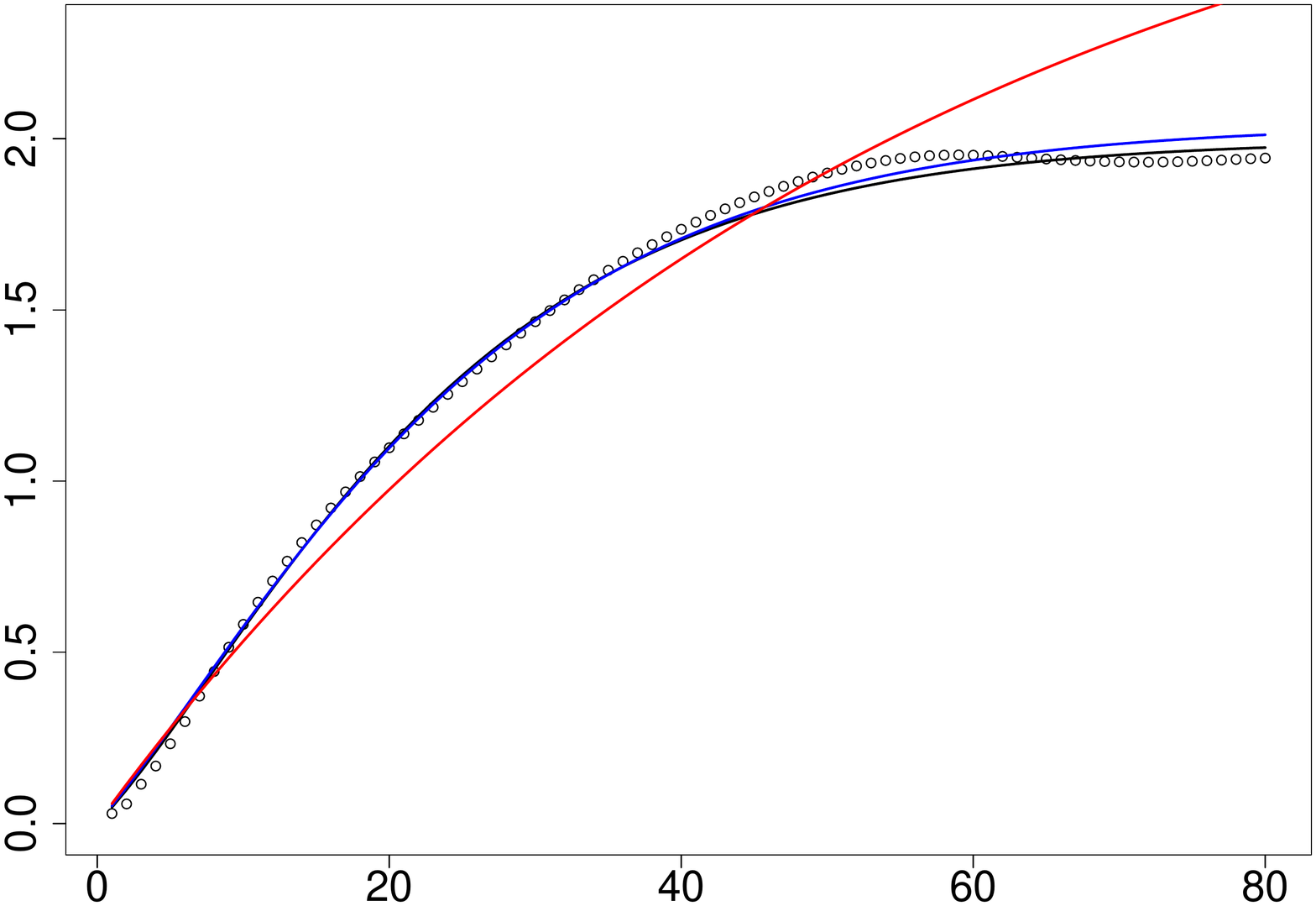}    
  \end{minipage}
  \begin{minipage}[b]{0.5\linewidth}
    \centering
    \includegraphics[width=\linewidth]{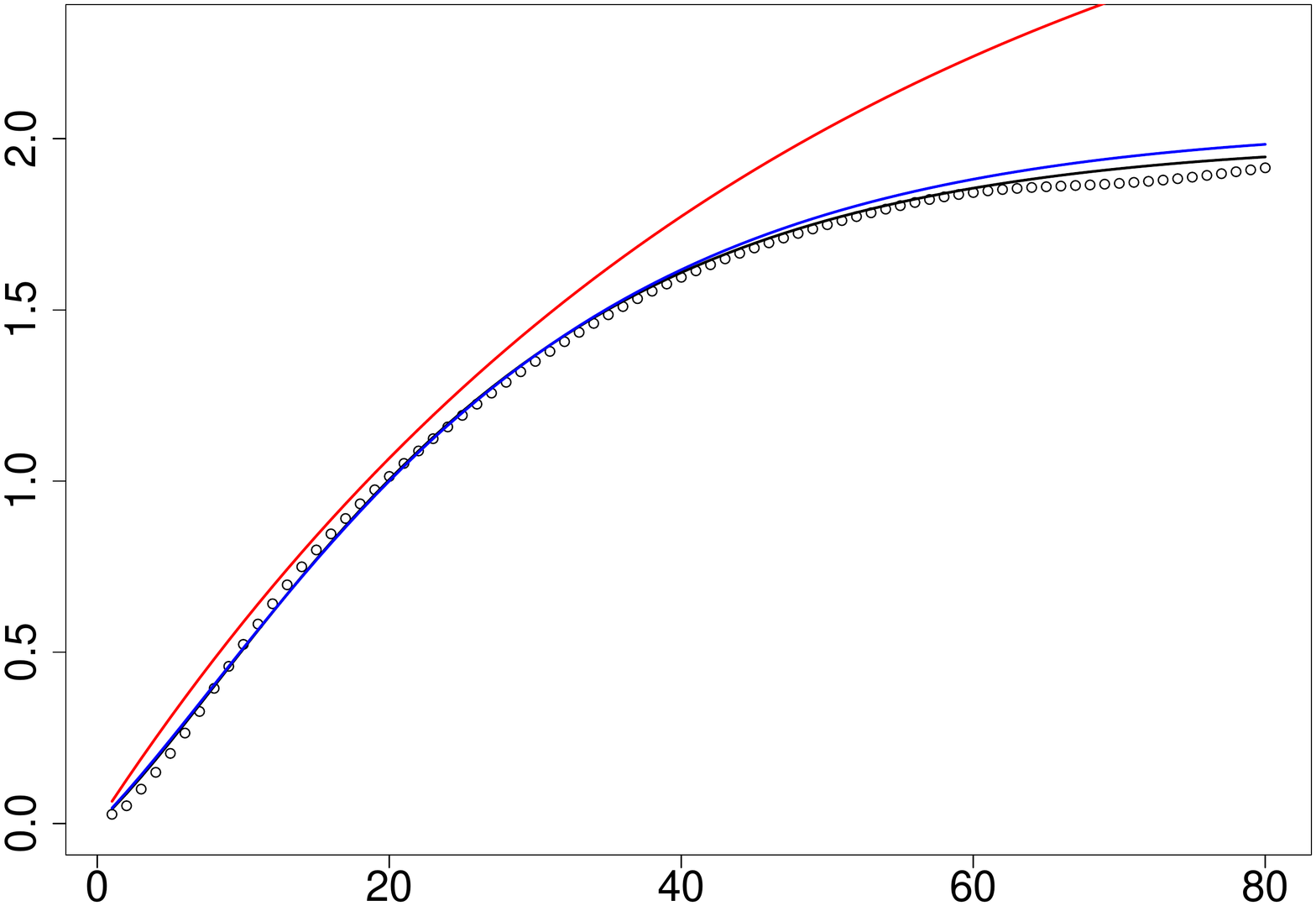}
  \end{minipage}
  \caption{Empirical variogram ordinates of the CMB data on the horizontal axis (left) and vertical axis (right) together with fitted variogram curves of the CAR$(1)$ model (in red), the CAR$(2)$ model (in blue) and the CARMA$(2,1)$ model (in black).} \label{variocompare}
\end{figure}

We apply our theory to cosmic microwave background (CMB) data from the Planck mission of the European Space Agency. The 2018 data release can be downloaded publicly from the Planck Legacy Archive \url{https://pla.esac.esa.int}. The CMB maps on this website cover the full sky and have been produced using four different methods. We choose the data set created by the SMICA method and refer to \cite{Akrami18} for more information. We take data points between $50^\circ$ and $70^\circ$ longitude and $10^\circ$ and $30^\circ$ latitude, the unit is given in Kelvin. We save the data with mean $-8.7316\times10^{-6}$ and standard deviation $9.6049\times10^{-5}$ into an $N\times N$-matrix with $N=1000$, and plot column-wise and row-wise means.
\begin{figure}[ht]
  \begin{minipage}[b]{0.5\linewidth}
    \centering
    \includegraphics[width=\linewidth]{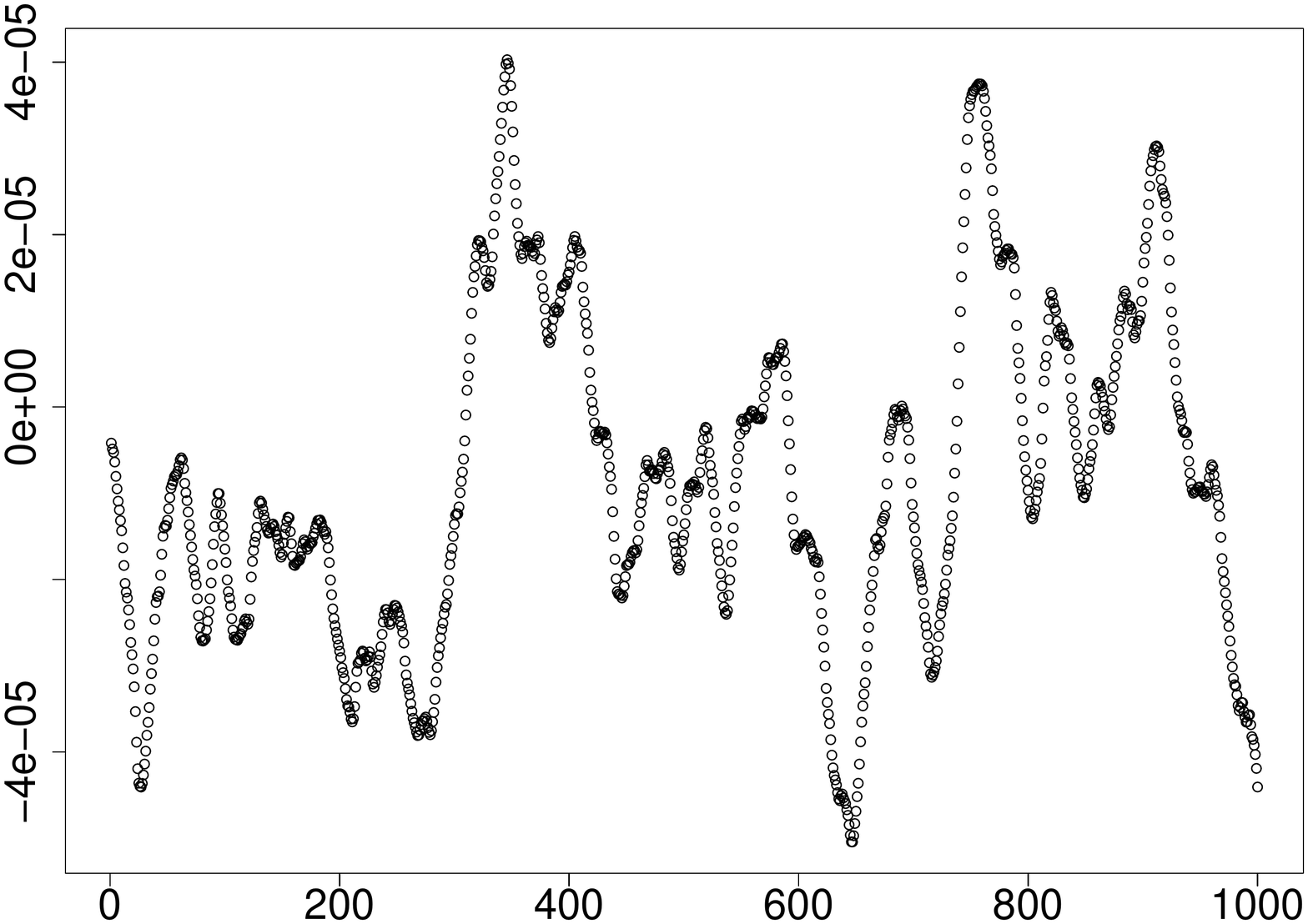}    
  \end{minipage}
  \begin{minipage}[b]{0.5\linewidth}
    \centering
    \includegraphics[width=\linewidth]{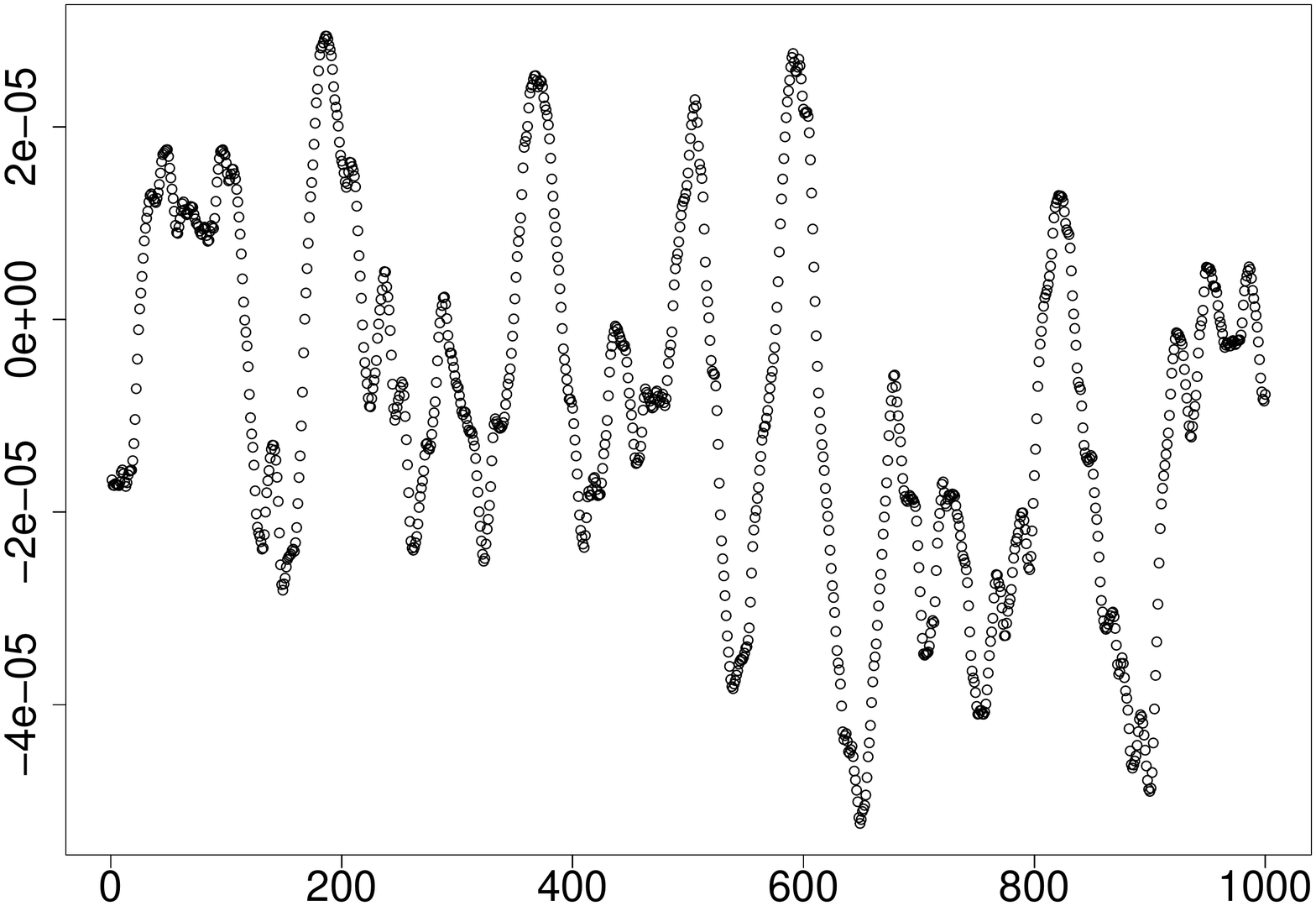}
  \end{minipage}
  \caption{Column-wise (left) and row-wise (right) means of the CMB data.}
  \label{colrowMeans}
\end{figure}
Since we do not find any deterministic trend or seasonal component in Figure~\ref{colrowMeans}, we may assume that the data is stationary. This is in line with standard assumptions for the CMB (see e.g. Section~2.1.1 of \citet{Giovannini08}).
We perform a normalization of the data to have mean zero and variance one, and plot the data's empirical density against the standard normal density.
\begin{figure}[ht]
  \centering
  \includegraphics[width=0.6\textwidth]{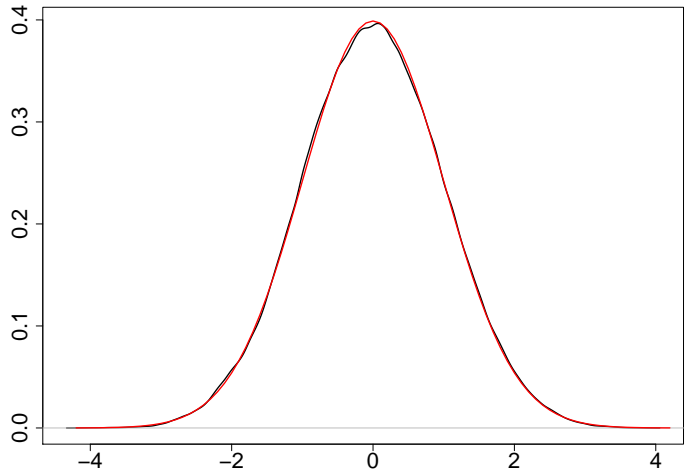}
  \caption{Empirical density of normalized CMB data (in black) and density of the standard normal distribution (in red).}
  \label{CMBdensity}
\end{figure}
An inspection of  Figure~\ref{CMBdensity} reveals that the marginal distribution of the CMB data is Gaussian. Hence, we may also assume that the Lévy basis $\La$ is Gaussian.
We proceed as in the previous section and compute the empirical variogram at $100$ different lags on the principal axes, namely $\{ \psi_{N}^*(j\Del\be_i)\colon i=1,2;j=1,...,50 \}$ with $\Del=0.04$. Assuming that the Lévy basis $\La$ has variance one, we estimate the parameters of CAR$(1)$, CAR$(2)$ and CARMA$(2,1)$ random fields with the WLS estimator in Equation~\eqref{case1}.
For the CAR$(1)$ model we obtain:
\begin{equation*}
b_0^*=1.2268,\, \la_{11}^*=-0.4622,\, \la_{21}^*=-0.5159.
\end{equation*}
For the CAR$(2)$ model we obtain:
\begin{equation*}
b_0^*=4.9991,\, \la_{11}^*=-1.7963,\, \la_{12}^*=-1.7969,\, \la_{21}^*=-1.2859,\, \la_{22}^*=-2.2212.
\end{equation*}
For the CARMA$(2,1)$ model we obtain:
\begin{equation*}
b_0^*=4.8940,\, b_1^*=-1.1432,\, \la_{11}^*=-1.7776,\, \la_{12}^*=-2.0948,\, \la_{21}^*=-1.3057,\, \la_{22}^*=-2.5142.
\end{equation*}
Figure~\ref{variocompare} depicts the estimated variogram of the CMB data along with fitted variogram curves of our three models. Recall that $b_0$ is not, but plays the role of the white noise standard deviation in our parametrization (see the first paragraph of Section~\ref{sectionestimation}).
The weighted sum of squares (WSS) values
\begin{equation*}
\text{WSS}=\sum_{\substack{j=1,...,50\\i=1,2}} w_j\left(\psi_{N}^*(j\Del\be_i)-\psi_{\theta^*}(j\Del\be_i)\right)^2
\end{equation*}
are $7.6132\times10^{-2}$ for CAR$(1)$, $2.5769\times10^{-2}$ for CAR$(2)$ and $2.0113\times10^{-2}$ for CARMA$(2,1)$. For model selection, we compute the Akaike information criterion (AIC)
\begin{equation*}
\text{AIC} = 2P + K\log(WSS/K),
\end{equation*}
where $P$ is the number of model parameters and $K$ the number of lags used to calculate the WSS. The AIC values are $-712.0453$ for CAR$(1)$, $-816.3761$ for CAR$(2)$ and $-839.1583$ for CARMA$(2,1)$.
\begin{table}[!htbp]
\centering
\begin{tabular}{rrrrr}
\hline
Model  & WSS & P  & K & AIC\\ 
\hline
CAR(1) & $7.6132\times10^{-2}$ & 3 & 100 & $-712.0453$ \\ 
CAR(2) & $2.5769\times10^{-2}$ & 5 & 100 & $-816.3761$ \\ 
CARMA(2,1) & $2.0113\times10^{-2}$ & 6 & 100 & $-839.1583$ \\ 
\hline
\end{tabular}
\caption{Weighted sum of squares (WSS), number of parameters (P), number of lags used (K) and Akaike Information Criterion (AIC) for the parameter estimation in Section~\ref{sectionCMB}.\label{table1}}
\end{table}

\pagebreak 

\begin{figure}[!htbp]
  \begin{minipage}[b]{0.5\linewidth}
    \centering
    \includegraphics[width=\linewidth]{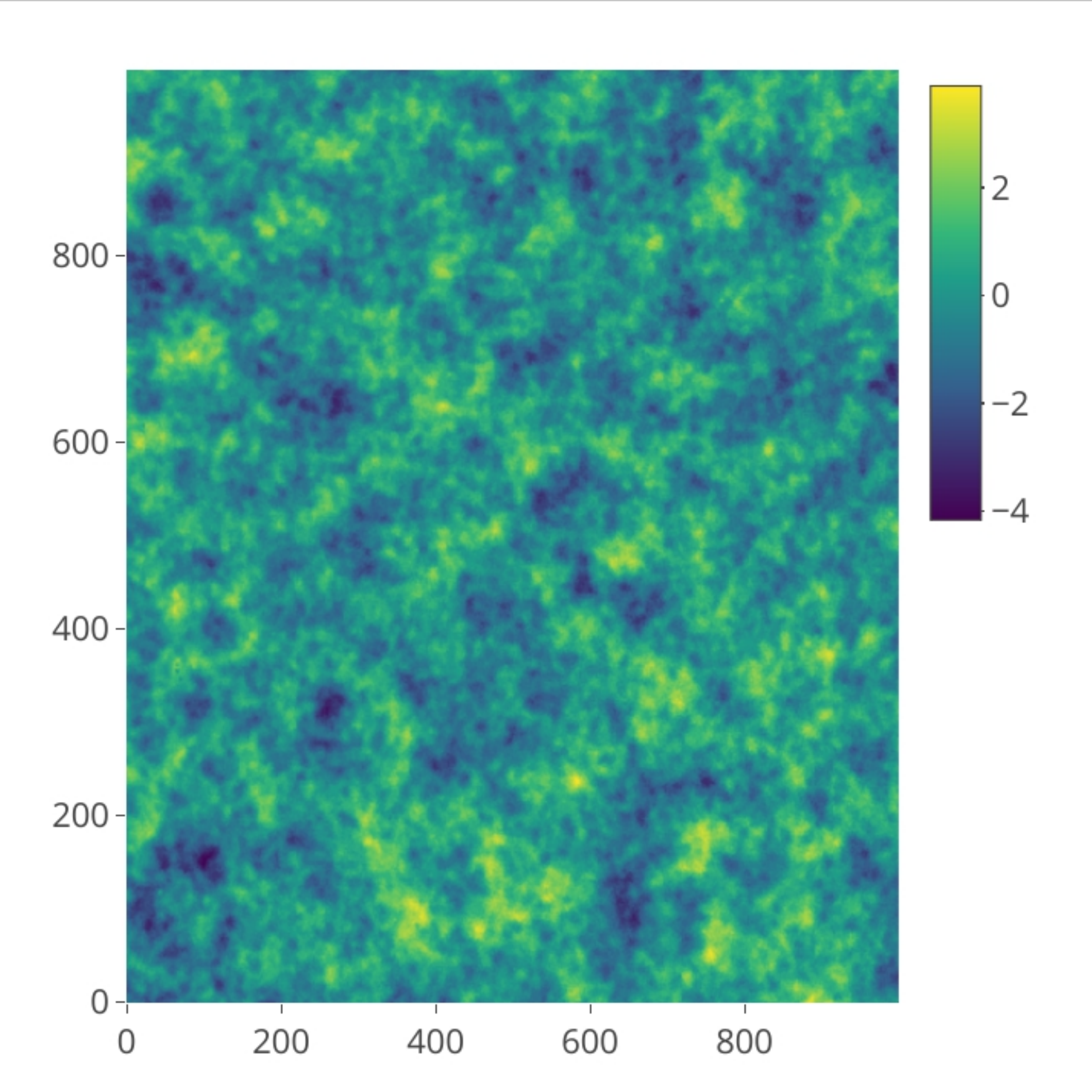}
    \caption{Normalized CMB data}
    \label{CMBdata.pdf}
  \end{minipage}
  \begin{minipage}[b]{0.5\linewidth}
    \centering
    \includegraphics[width=\linewidth]{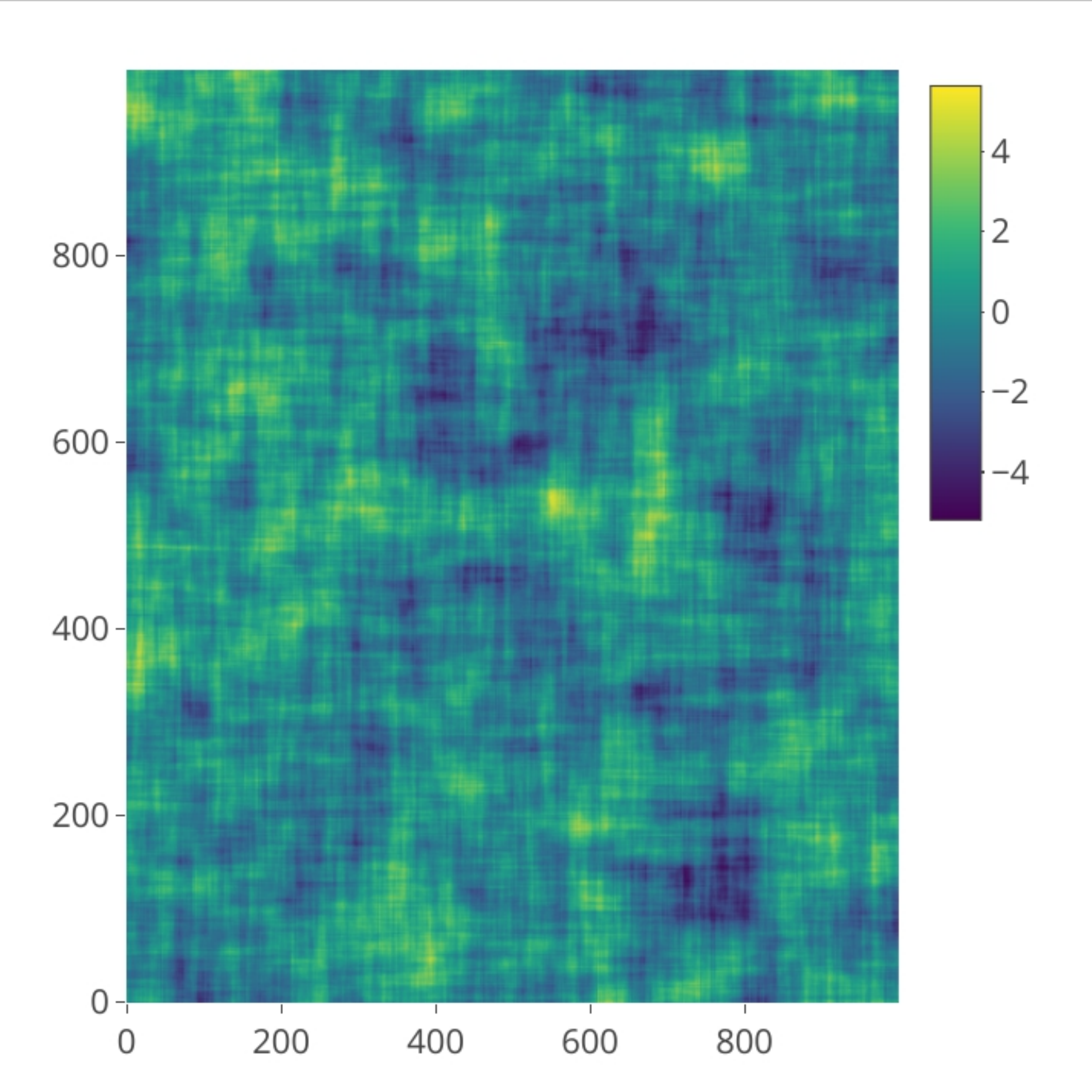}
    \caption{Simulated CAR$(1)$ random field}
    \label{CAR(1).pdf}
  \end{minipage}
  \begin{minipage}[b]{0.5\linewidth}
      \centering
      \includegraphics[width=\linewidth]{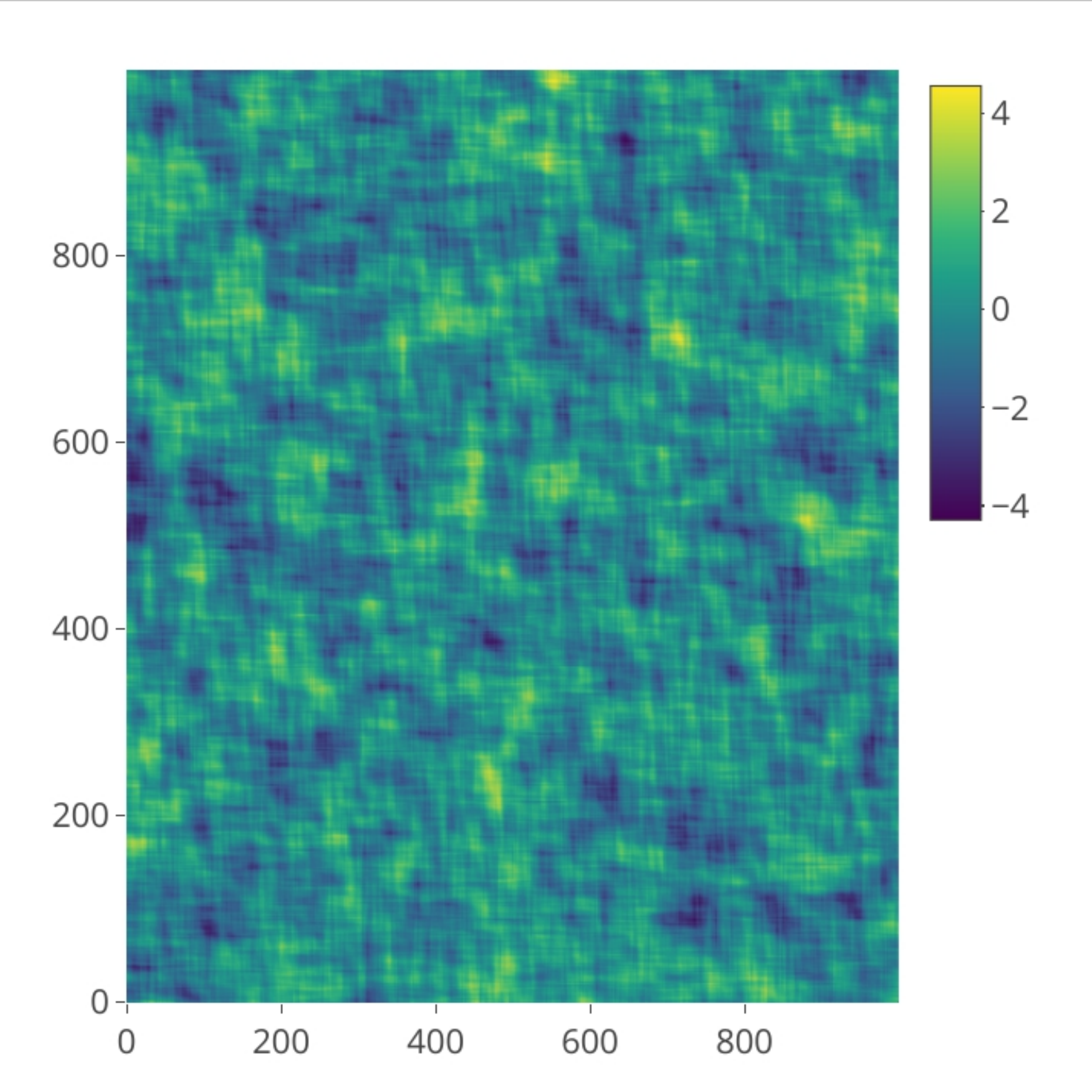}
      \caption{Simulated CAR$(2)$ random field}
      \label{CAR(2).pdf}
    \end{minipage}
    \begin{minipage}[b]{0.5\linewidth}
      \centering
      \includegraphics[width=\linewidth]{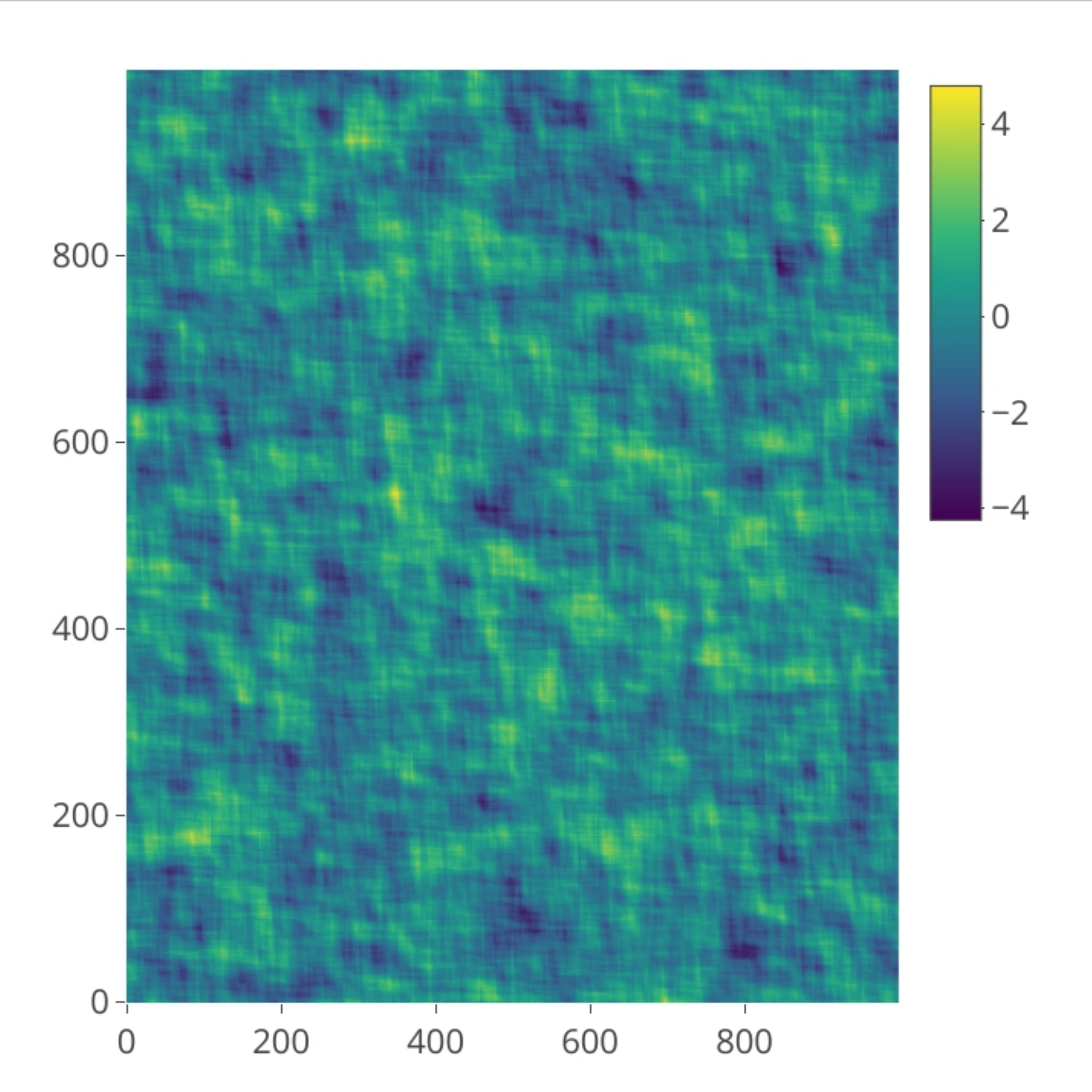}
      \caption{Simulated CARMA$(2,1)$ random field}
      \label{CARMA(2,1).pdf}
    \end{minipage}
\end{figure}

These numbers are summarized in Table~\ref{table1} and suggest that the CARMA$(2,1)$ model is optimal compared to the CAR$(1)$ and CAR$(2)$ models. For a visual comparison we plot the heat map of the original CMB data together with heat maps of simulated fields in Figures~\ref{CMBdata.pdf} to \ref{CARMA(2,1).pdf}. Although we cannot draw any conclusions from a single sample path, it is possible to observe some features of the fitted models. All three models exhibit clusters of high and low values similarly to the original data. However, the cluster sizes of the CAR$(1)$ random field are larger than those of the CMB data, whereas the CAR$(2)$ and CARMA$(2,1)$ models display a better visual fit. Another common feature are horizontal and vertical lines, which is most visible in Figure~\ref{CAR(1).pdf}. These lines are the consequences of the non-smoothness of the kernel function in Equation~\eqref{kernel}. One therefore can argue that the fitted CARMA random fields represent linear approximations to the spatial dependence structures of the cosmic microwave background.

\begin{table}[ht]
\centering
\begin{tabular}{rrrrrr}
  \hline
  & True Value & Mean & Bias & Std & RMSE \\ 
  \hline
$b_0$ & 4.8940 & 4.7882 & -0.1058 & 0.5124 & 0.5227 \\ 
$b_1$ & -1.1432 & -1.2784 & -0.1352 & 0.3962 & 0.4183 \\ 
$\la_{11}$ & -1.7776 & -1.6283 & 0.1494 & 0.2377 & 0.2806 \\ 
$\la_{12}$ & -2.0948 & -2.3193 & -0.2246 & 0.4183 & 0.4744 \\ 
$\la_{21}$ & -1.3057 & -1.3136 & -0.0079 & 0.2323 & 0.2322 \\ 
$\la_{22}$ & -2.5142 & -2.5231 & -0.0089 & 0.4048 & 0.4045 \\ 
   \hline
\end{tabular}
\caption{Parameter estimation results for CARMA$(2,1)$ on $\bbr^2$ with $K=100$ lags, quadratically decreasing weights as in \eqref{case1} and Gaussian basis $\La$.\label{summary1}}
\begin{tabular}{rrrrrr}
  \hline
 & True Value & Mean & Bias & Std & RMSE \\ 
  \hline
$b_0$   & 4.8940 & 4.6929 & -0.2010 & 0.4597 & 0.5013 \\ 
$b_1$   & -1.1432 & -1.2252 & -0.0820 & 0.3515 & 0.3606 \\ 
$\la_{11}$   & -1.7776 & -1.6335 & 0.1442 & 0.2005 & 0.2468 \\ 
$\la_{12}$   & -2.0948 & -2.2117 & -0.1169 & 0.3246 & 0.3447 \\ 
$\la_{21}$   & -1.3057 & -1.2947 & 0.0110 & 0.2136 & 0.2137 \\ 
$\la_{22}$   & -2.5142 & -2.4636 & 0.0506 & 0.3065 & 0.3104 \\ 
   \hline
\end{tabular}
\caption{Parameter estimation results for CARMA$(2,1)$ on $\bbr^2$ with $K=50$ lags, quadratically decreasing weights as in \eqref{case2} and Gaussian basis $\La$.\label{summary2}}
\begin{tabular}{rrrrrr}
  \hline
 & True Value & Mean & Bias & Std & RMSE \\ 
  \hline
$b_0$   & 4.8940 & 4.8329 & -0.0610 & 0.5668 & 0.5695 \\ 
$b_1$   & -1.1432 & -1.2708 & -0.1275 & 0.4250 & 0.4433 \\ 
$\la_{11}$   & -1.7776 & -1.6234 & 0.1542 & 0.2473 & 0.2912 \\ 
$\la_{12}$   & -2.0948 & -2.3569 & -0.2622 & 0.4754 & 0.5425 \\ 
$\la_{21}$   & -1.3057 & -1.3182 & -0.0125 & 0.2448 & 0.2449 \\ 
$\la_{22}$   & -2.5142 & -2.5392 & -0.0250 & 0.4348 & 0.4351 \\ 
   \hline
\end{tabular}
\caption{Parameter estimation results for CARMA$(2,1)$ on $\bbr^2$ with $K=100$ lags, exponentially decreasing weights as in \eqref{case3} and Gaussian basis $\La$.\label{summary3}}
\begin{tabular}{rrrrrr}
  \hline
 & True Value & Mean & Bias & Std & RMSE \\ 
  \hline
$b_0$   & 4.8940 & 4.7525 & -0.1414 & 0.5267 & 0.5448 \\ 
$b_1$   & -1.1432 & -1.1995 & -0.0563 & 0.3879 & 0.3915 \\ 
$\la_{11}$   & -1.7776 & -1.5908 & 0.1868 & 0.2375 & 0.3020 \\ 
$\la_{12}$   & -2.0948 & -2.2988 & -0.2040 & 0.4240 & 0.4701 \\ 
$\la_{21}$   & -1.3057 & -1.2765 & 0.0292 & 0.2299 & 0.2315 \\ 
$\la_{22}$   & -2.5142 & -2.5196 & -0.0054 & 0.3786 & 0.3783 \\ 
   \hline
\end{tabular}
\caption{Parameter estimation results for CARMA$(2,1)$ on $\bbr^2$ with $K=50$ lags, exponentially decreasing weights as in \eqref{case4} and Gaussian basis $\La$.\label{summary4}}
\end{table}

\begin{table}[ht]
\centering
\begin{tabular}{rrrrrr}
  \hline
 & True Value & Mean & Bias & Std & RMSE \\ 
  \hline
$b_0$   & 4.8940 & 4.8407 & -0.0533 & 0.5126 & 0.5148 \\ 
$b_1$   & -1.1432 & -1.2383 & -0.0951 & 0.4181 & 0.4283 \\ 
$\la_{11}$   & -1.7776 & -1.6453 & 0.1323 & 0.2202 & 0.2567 \\ 
$\la_{12}$   & -2.0948 & -2.3054 & -0.2106 & 0.3828 & 0.4366 \\ 
$\la_{21}$   & -1.3057 & -1.3149 & -0.0092 & 0.2269 & 0.2269 \\ 
$\la_{22}$   & -2.5142 & -2.5319 & -0.0177 & 0.3716 & 0.3717 \\ 
   \hline
\end{tabular}
\caption{Parameter estimation results for CARMA$(2,1)$ on $\bbr^2$ with $K=100$ lags, quadratically decreasing weights as in \eqref{case1} and variance gamma basis $\La$.\label{summary5}}
\begin{tabular}{rrrrrr}
  \hline
 & True Value & Mean & Bias & Std & RMSE \\ 
  \hline
$b_0$ & 4.8940 & 4.7474 & -0.1466 & 0.4635 & 0.4857 \\ 
$b_1$ & -1.1432 & -1.2035 & -0.0603 & 0.3578 & 0.3625 \\ 
$\la_{11}$ & -1.7776 & -1.6325 & 0.1451 & 0.2034 & 0.2497 \\ 
$\la_{12}$ & -2.0948 & -2.2322 & -0.1375 & 0.3246 & 0.3522 \\ 
$\la_{21}$ & -1.3057 & -1.2866 & 0.0191 & 0.2137 & 0.2144 \\ 
$\la_{22}$ & -2.5142 & -2.4994 & 0.0148 & 0.3154 & 0.3154 \\ 
   \hline
\end{tabular}
\caption{Parameter estimation results for CARMA$(2,1)$ on $\bbr^2$ with $K=50$ lags, quadratically decreasing weights as in \eqref{case2} and variance gamma basis $\La$.\label{summary6}}
\begin{tabular}{rrrrrr}
  \hline
 & True Value & Mean & Bias & Std & RMSE \\ 
  \hline
$b_0$ & 4.8940 & 4.8696 & -0.0243 & 0.5671 & 0.5671 \\ 
$b_1$ & -1.1432 & -1.2780 & -0.1348 & 0.4090 & 0.4302 \\ 
$\la_{11}$ & -1.7776 & -1.6266 & 0.1511 & 0.2543 & 0.2956 \\ 
$\la_{12}$ & -2.0948 & -2.3773 & -0.2825 & 0.4731 & 0.5506 \\ 
$\la_{21}$ & -1.3057 & -1.3086 & -0.0029 & 0.2402 & 0.2400 \\ 
$\la_{22}$ & -2.5142 & -2.5760 & -0.0618 & 0.4262 & 0.4303 \\ 
   \hline
\end{tabular}
\caption{Parameter estimation results for CARMA$(2,1)$ on $\bbr^2$ with $K=100$ lags, exponentially decreasing weights as in \eqref{case3} and variance gamma basis $\La$.\label{summary7}}
\begin{tabular}{rrrrrr}
  \hline
 & True Value & Mean & Bias & Std & RMSE \\ 
  \hline
$b_0$ & 4.8940 & 4.7422 & -0.1518 & 0.4806 & 0.5035 \\ 
$b_1$ & -1.1432 & -1.2331 & -0.0898 & 0.3942 & 0.4039 \\ 
$\la_{11}$ & -1.7776 & -1.6102 & 0.1674 & 0.2193 & 0.2758 \\ 
$\la_{12}$ & -2.0948 & -2.2791 & -0.1843 & 0.3626 & 0.4065 \\ 
$\la_{21}$ & -1.3057 & -1.2816 & 0.0241 & 0.2272 & 0.2283 \\ 
$\la_{22}$ & -2.5142 & -2.5225 & -0.0083 & 0.3532 & 0.3530 \\ 
   \hline
\end{tabular}
\caption{Parameter estimation results for CARMA$(2,1)$ on $\bbr^2$ with $K=50$ lags, exponentially decreasing weights as in \eqref{case4} and variance gamma basis $\La$.\label{summary8}}
\end{table}

\subsection*{Acknowledgement}

We cordially thank Thiago do Rego Sousa for helpful discussions. The second author acknowledges support from the graduate program TopMath at the Technical University of Munich and the Studienstiftung des deutschen Volkes.



\end{document}